\documentclass{amsart}
\usepackage{amsfonts}

\newtheorem{theorem}{Theorem}
\theoremstyle{plain}

\newtheorem{definition}{Definition}

\newtheorem{lemma}{Lemma}

\numberwithin{equation}{section}
\input{tcilatex}

\begin{document}
\title[Everywhere H\"{o}lder continuity]{On the local everywhere H\"{o}lder
continuity for weak solutions of a class of not convex vectorial problems of
the Calculus of Variations}
\author{Tiziano Granucci}
\address{MIUR, Istituto Leonardo da Vinci, via del Terzolle 91, 50127,
Firenze, Italy}
\email{tizianogranucci@libero.it}
\urladdr{https://www.tizianogranucci.com/}
\thanks{I thank my family Elisa Cirri, Caterina Granucci, Delia Granucci for
their support. I also thank my friends Monia Randolfi and Massimo Masi for
the many discussions and for the many advice.}
\date{27/06/2022}
\dedicatory{Dedicated to my family, Elisa Cirri, Caterina Granucci and Delia
Granucci}
\thanks{This paper is in final form and no version of it will be submitted
for publication elsewhere.}

\begin{abstract}
In this paper we study the regularity of the local minima of the following
integral functional%
\begin{equation}
J\left( u,\Omega \right) =\int\limits_{\Omega }\sum\limits_{\alpha
=1}^{n}\left\vert \nabla u^{\alpha }\left( x\right) \right\vert ^{p}+G\left(
x,u\left( x\right) ,\nabla u\left( x\right) \right) \,dx
\end{equation}%
where $\Omega $ is a open subset of $%
%TCIMACRO{\U{211d} }%
%BeginExpansion
\mathbb{R}
%EndExpansion
^{n}$ and $u\in W^{1,p}\left( \Omega ,%
%TCIMACRO{\U{211d} }%
%BeginExpansion
\mathbb{R}
%EndExpansion
^{m}\right) $ with $n\geq 2$, $m\geq 1$ and $1<p<n$. In particular, not
convexity (quasi-convexity, policonvexity or rank one convexity) hypothesis
will be made on the density $G$, neither structure hypothesis nor radial nor
diagonal.
\end{abstract}

\maketitle

\section{Introduction}

In this paper we study the regularity of the local minima of the following
integral functional%
\begin{equation}
J\left( u,\Omega \right) =\int\limits_{\Omega }\sum\limits_{\alpha
=1}^{n}\left\vert \nabla u^{\alpha }\left( x\right) \right\vert ^{p}+G\left(
x,u\left( x\right) ,\nabla u\left( x\right) \right) \,dx  \label{1.1}
\end{equation}%
where $\Omega $ is a open subset of $%
%TCIMACRO{\U{211d} }%
%BeginExpansion
\mathbb{R}
%EndExpansion
^{n}$ and $u\in W^{1,p}\left( \Omega ,%
%TCIMACRO{\U{211d} }%
%BeginExpansion
\mathbb{R}
%EndExpansion
^{m}\right) $ with $n\geq 2$, $m\geq 1$ and $1<p<n$.

Moreover the following hypotheses hold

\begin{description}
\item[H.1.1] \bigskip $G:\Omega \times 
%TCIMACRO{\U{211d} }%
%BeginExpansion
\mathbb{R}
%EndExpansion
^{m}\times 
%TCIMACRO{\U{211d} }%
%BeginExpansion
\mathbb{R}
%EndExpansion
^{n\times m}\rightarrow 
%TCIMACRO{\U{211d} }%
%BeginExpansion
\mathbb{R}
%EndExpansion
$ is a Caratheodory function such that%
\begin{equation*}
\left\vert G\left( x,s_{1},\xi _{1}\right) -G\left( x,s_{2},\xi _{2}\right)
\right\vert \leq a\left( x\right) \left\vert s_{1}-s_{2}\right\vert ^{\alpha
}\left( \left\vert \xi _{1}\right\vert +\left\vert \xi _{2}\right\vert
+1\right) ^{q_{1}}+b\left( x\right) \left( \left\vert s_{1}\right\vert
+\left\vert s_{2}\right\vert +1\right) ^{q_{2}}\left\vert \xi _{1}-\xi
_{2}\right\vert ^{\beta }
\end{equation*}%
for $\mathcal{L}^{n}$ almost every $x\in \Omega $ and for every $%
s_{1},s_{2}\in 
%TCIMACRO{\U{211d} }%
%BeginExpansion
\mathbb{R}
%EndExpansion
^{m}$ and $\xi _{1},\xi _{2}\in 
%TCIMACRO{\U{211d} }%
%BeginExpansion
\mathbb{R}
%EndExpansion
^{n\times m}$, where $0<\alpha <\min \left\{ 1,p-q_{1}\right\} $, $q_{1}<p$, 
$0<\beta <\min \left\{ 1,p-q_{2}\right\} $, $q_{2}<p$, $a\in L_{loc}^{\sigma
_{1}}\left( \Omega \right) $ is a not negative function, $b\in
L_{loc}^{\sigma _{2}}\left( \Omega \right) $ is a not negative function, $%
\sigma _{1}>\frac{pp^{\ast }}{\left( p^{\ast }-\alpha \right) p-q_{1}p^{\ast
}}$, $\sigma _{2}>\frac{pp^{\ast }}{\left( p-\beta \right) p^{\ast }-q_{2}p}$%
, \ 
\begin{equation*}
\frac{q_{1}}{p}+\frac{\alpha }{p^{\ast }}+\frac{1}{\sigma _{1}}<\frac{p}{n}
\end{equation*}%
and%
\begin{equation*}
\frac{\beta }{p}+\frac{q_{2}}{p^{\ast }}+\frac{1}{\sigma _{2}}<\frac{p}{n}
\end{equation*}
\end{description}

The main result of this article is the following regularity theorem:

\begin{theorem}
If $u\in W^{1,p}\left( \Omega ,%
%TCIMACRO{\U{211d} }%
%BeginExpansion
\mathbb{R}
%EndExpansion
^{m}\right) $ is a minimizer of (1.1) and H.1.1 holds then $u\in
C_{loc}^{0,\delta }\left( \Omega ,%
%TCIMACRO{\U{211d} }%
%BeginExpansion
\mathbb{R}
%EndExpansion
^{m}\right) $.
\end{theorem}

\bigskip Theorem 1 is interesting for a few reasons. We know that in the
vector case there are many counter examples, refer to [14, 19, 21], and in
general the minima are not everywhere regular, refeer to [16, 37].
Furthermore, starting from the end of the 1970s, using suitable hypotheses
of convexity and regularity on the density $\Phi $ for the minima of
functionals of the type 
\begin{equation*}
\int\limits_{\Omega }\Phi \left( \left\vert \nabla u\left( x\right)
\right\vert \right) \,dx
\end{equation*}%
regularity theorems have been proved, refer to [1, 2, 4-7, 15, 17, 18, 20,
21, 34, 42-44]. Theorem 1 has no hypothesis either of structure or form, or
of regularity or convesity on the density $G$. Finally, the proof of Theorem
1 is particularly simple, in fact the previous Theorem 1 derives from the
following Cacciopoli inequalities using the techniques introduced by E. De
Giorgi in [13].

\begin{theorem}
If $u\in W^{1,p}\left( \Omega ,%
%TCIMACRO{\U{211d} }%
%BeginExpansion
\mathbb{R}
%EndExpansion
^{m}\right) $ is a minimizer of () and H.1.1 holds then, for every $\Sigma
\subset \Omega $ compact, two positive constants $C_{Cac,1}$, $C_{Cac,2}$
(dipendenti solo da $\Sigma $, $p$ e $n$) and a radius $R_{0}>0$ exist such
that for every $0<\varrho <R<R_{0}$ for every $x_{0}\in \Sigma $ \ and for
every $k\in 
%TCIMACRO{\U{211d} }%
%BeginExpansion
\mathbb{R}
%EndExpansion
$ it follows%
\begin{equation*}
\int\limits_{A_{k,\varrho }^{\alpha }}\left\vert \nabla u^{\alpha
}\right\vert ^{p}\,dx\leq \frac{C_{Cac,1}}{\left( R-\varrho \right) ^{p}}%
\int\limits_{A_{k,R}^{\alpha }}\left( u^{\alpha }-k\right)
^{p}\,dx+C_{Cac,2}\left\vert A_{k,s}^{\alpha }\right\vert ^{1-\frac{p}{N}%
+\epsilon }
\end{equation*}%
and%
\begin{equation*}
\int\limits_{B_{k,\varrho }^{\alpha }}\left\vert \nabla u^{\alpha
}\right\vert ^{p}\,dx\leq \frac{C_{Cac,1}}{\left( R-\varrho \right) ^{p}}%
\int\limits_{B_{k,R}^{\alpha }}\left( k-u^{\alpha }\right)
^{p}\,dx+C_{Cac,2}\left\vert B_{k,s}^{\alpha }\right\vert ^{1-\frac{p}{N}%
+\epsilon }
\end{equation*}%
where $A_{k,s}^{\alpha }=\left\{ u^{\alpha }>k\right\} \cap B_{s}\left(
x_{0}\right) $ and $B_{k,s}^{\alpha }=\left\{ u^{\alpha }<k\right\} \cap
B_{s}\left( x_{0}\right) $ with $\alpha =1,...,m$.
\end{theorem}

Our results can therefore be framed within a vast area of {}{}research
called everywhere regulairy that was born with the fundamental works of
Uhlenbeck [44], Tolksdorf [42, 43] and Acerbi - Fusco [1]. The literature in
this area is very wide, in the bibliography we report only some articles,
refer to[1,2, 4-7, 15, 17-21, 34-36, 42-44]. Theorem 1 is part of the study
of the regularity of the minima of functionals of the following type%
\begin{equation*}
\int\limits_{B_{R}\left( x_{0}\right) }F\left( x,u,\nabla u\right) \,dx
\end{equation*}%
in this case, as far as the limited knowledge of the author is concerned,
there are few results, in particular we refer to [8-11, 25, 29-33]. The
previous theorem 1 can therefore be included in a series of results obtained
in recent years [8-11, 25, 29-33]. In [8] Cupini, Focardi, Leonetti and
Mascolo introduced the following class of vectorial functionals%
\begin{equation}
\int\limits_{\Omega }f\left( x,\nabla u\right) \,dx  \label{1.3}
\end{equation}%
where $\Omega \subset 
%TCIMACRO{\U{211d} }%
%BeginExpansion
\mathbb{R}
%EndExpansion
^{n},u:\Omega \rightarrow 
%TCIMACRO{\U{211d} }%
%BeginExpansion
\mathbb{R}
%EndExpansion
^{m},n>1,m\geq 1$ and%
\begin{equation*}
f\left( x,\nabla u\right) =\sum\limits_{\alpha =1}^{m}F_{\alpha }\left(
x,\nabla u^{\alpha }\right) +G\left( x,\nabla u\right)
\end{equation*}%
where $F\alpha \times \ R^{n\times m}\rightarrow R$ is a Carath\'{e}odory
function satisfying the following standard growth condition%
\begin{equation*}
k_{1}\left\vert \xi ^{\alpha }\right\vert ^{p}-a\left( x\right) \leq
F_{\alpha }\left( x,\xi ^{\alpha }\right) \leq k_{2}\left\vert \xi ^{\alpha
}\right\vert ^{p}+a\left( x\right)
\end{equation*}%
for every $\xi ^{\alpha }\in 
%TCIMACRO{\U{211d} }%
%BeginExpansion
\mathbb{R}
%EndExpansion
^{n}$ and for almost every $x\in \Omega $ ,where $k_{1}$ and $k_{2}$ are two
real positive constants, $p>1$ and $a\in L_{loc}^{\sigma }\left( \Omega
\right) $ is a non negative function. In [8],Cupini, Focardi, Leonetti and
Mascolo analyze two different types of hypotheses on the $G$ function. They
started by assuming that $G:\Omega \times \ R^{n\times m}\rightarrow R$ is a
Carath\'{e}odory rank one convex function satisfying the following growth
condition%
\begin{equation*}
|G(x,\xi )|\leq k_{3}|\xi |^{q}+b(x)
\end{equation*}%
for every $\xi \in 
%TCIMACRO{\U{211d} }%
%BeginExpansion
\mathbb{R}
%EndExpansion
^{n\times m}$, for almost every $x\in \Omega $, here $k_{3}$ is a real
positive constant, $1\leq q<p$ and $b\in L_{loc}^{\sigma }\left( \Omega
\right) $ is a nonnegative function. Moreover Cupini, Focardi, Leonetti and
Mascolo in [8] study the case where $n\geq m\geq 3$, and $G:\Omega \times \ 
%TCIMACRO{\U{211d} }%
%BeginExpansion
\mathbb{R}
%EndExpansion
^{n\times m}\rightarrow R$ is a Carath\'{e}odory function defined as%
\begin{equation*}
G(x,\xi )=\sum\limits_{\alpha =1}^{m}G_{\alpha }\left( x,\left( adj_{m-1}\xi
\right) ^{\alpha }\right)
\end{equation*}%
here $G_{\alpha }:\Omega \times \ 
%TCIMACRO{\U{211d} }%
%BeginExpansion
\mathbb{R}
%EndExpansion
^{\frac{m!}{n!\left( n-m\right) !}}\rightarrow 
%TCIMACRO{\U{211d} }%
%BeginExpansion
\mathbb{R}
%EndExpansion
$ is a Carath\'{e}odory convex function satisfying the following growth
conditions%
\begin{equation*}
0\leq G_{\alpha }\left( x,\left( adj_{m-1}\xi \right) ^{\alpha }\right) \leq
k_{4}|\left( adj_{m-1}\xi \right) ^{\alpha }|^{r}+b(x)
\end{equation*}%
for every $\xi \in 
%TCIMACRO{\U{211d} }%
%BeginExpansion
\mathbb{R}
%EndExpansion
^{n\times m}$, for almost every $x\in \Omega $, here $k_{3}$ is a real
positive constant, $1\leq r<p$ and $b\in L_{loc}^{\sigma }\left( \Omega
\right) $ is a non negative function. In both cases, by imposing appropriate
hypotheses on the parameters $q$ and $r$ , Cupini, Focardi, Leonetti and
Mascolo proved that the local minimizers of the vectorial functional (\ref%
{1.3}) are locally h\"{o}lder continuous functions.

In [32] the author studied the regularity of the minima of the following
functional 
\begin{equation}
J\left( u,\Omega \right) =\int\limits_{\Omega }\sum\limits_{\alpha
=1}^{m}f_{\alpha }\left( x,u^{\alpha }\left( x\right) ,\nabla u^{\alpha
}\left( x\right) \right) +G\left( x,u\left( x\right) ,\nabla u\left(
x\right) \right) \,dx
\end{equation}%
where $\Omega $ is a open subset of $%
%TCIMACRO{\U{211d} }%
%BeginExpansion
\mathbb{R}
%EndExpansion
^{n}$ and $u\in W^{1,p}\left( \Omega ,%
%TCIMACRO{\U{211d} }%
%BeginExpansion
\mathbb{R}
%EndExpansion
^{m}\right) $ with $n\geq 2$, $m\geq 1$ and $1<p<n$ e supponendo che le
seguenti ipotesi valgano

\begin{description}
\item[H.2.1] For every $\alpha =1,...,m$ the function$\ f_{\alpha }:\Omega
\times 
%TCIMACRO{\U{211d} }%
%BeginExpansion
\mathbb{R}
%EndExpansion
\times 
%TCIMACRO{\U{211d} }%
%BeginExpansion
\mathbb{R}
%EndExpansion
^{n}\rightarrow 
%TCIMACRO{\U{211d} }%
%BeginExpansion
\mathbb{R}
%EndExpansion
$ is a Caratheodory function and the following growth conditions hold%
\begin{equation}
\left\vert \xi ^{\alpha }\right\vert ^{p}-b_{\alpha }\left( x\right)
\left\vert s\right\vert ^{\gamma _{\alpha }}-a_{\alpha }\left( x\right) \leq
f_{\alpha }\left( x,s,\xi ^{\alpha }\right) \leq L_{\alpha }\left(
\left\vert \xi ^{\alpha }\right\vert ^{p}+b_{\alpha }\left( x\right)
\left\vert s\right\vert ^{\gamma _{\alpha }}+a_{\alpha }\left( x\right)
\right)
\end{equation}%
for almost every $x\in \Omega $, for every $s\in 
%TCIMACRO{\U{211d} }%
%BeginExpansion
\mathbb{R}
%EndExpansion
$ and for every $\xi ^{\alpha }\in 
%TCIMACRO{\U{211d} }%
%BeginExpansion
\mathbb{R}
%EndExpansion
^{n}$ where $L_{\alpha }>1$, $1<p\leq \gamma _{\alpha }<p^{\ast }=\frac{np}{%
n-p}$, $b_{\alpha }$ and $a_{\alpha }$ are two not-negative function, $%
b_{\alpha }\in L_{loc}^{\sigma _{\alpha }}\left( \Omega \right) $ and $%
a_{\alpha }\in L_{loc}^{\kappa }\left( \Omega \right) $ with $\sigma
_{\alpha }=\frac{p^{\ast }}{p^{\ast }-\gamma _{\alpha }-\epsilon p^{\ast }}$%
, $\kappa =\frac{n}{p-\epsilon n}$ and $0<\epsilon <\frac{p}{n}$.

\item[H.2.2] $G:\Omega \times 
%TCIMACRO{\U{211d} }%
%BeginExpansion
\mathbb{R}
%EndExpansion
^{m}\times 
%TCIMACRO{\U{211d} }%
%BeginExpansion
\mathbb{R}
%EndExpansion
^{nm}\rightarrow 
%TCIMACRO{\U{211d} }%
%BeginExpansion
\mathbb{R}
%EndExpansion
$ is a Caratheodory function and the following growth conditions hold%
\begin{equation}
\left\vert G\left( x,u,\xi \right) \right\vert \leq C\left( \left\vert \xi
\right\vert ^{q}+\left\vert u\right\vert ^{q}+a\left( x\right) \right)
\end{equation}%
for almost every $x\in \Omega $, for every $u\in 
%TCIMACRO{\U{211d} }%
%BeginExpansion
\mathbb{R}
%EndExpansion
^{m}$ and for every $\xi \in 
%TCIMACRO{\U{211d} }%
%BeginExpansion
\mathbb{R}
%EndExpansion
^{nm}$ where $C>1$, $1\leq q<\frac{p^{2}}{n}<p$, $a$ is a not-negative
function and $a\in L_{loc}^{\kappa }\left( \Omega \right) $ with $\kappa =%
\frac{n}{p-\epsilon n}$ and $0<\epsilon <\frac{p}{n}$.

\item[H.2.3] The function $G\left( x,u,\cdot \right) $ is rank one convex
then%
\begin{equation*}
G\left( x,u,\lambda \xi ^{1}+\left( 1-\lambda \right) \xi ^{2}\right) \leq
\lambda G\left( x,u,\xi ^{1}\right) +\left( 1-\lambda \right) G\left(
x,u,\xi ^{2}\right)
\end{equation*}%
for a. e. $x\in \Omega $, for every $u\in 
%TCIMACRO{\U{211d} }%
%BeginExpansion
\mathbb{R}
%EndExpansion
^{m}$, for every $\lambda \in \left[ 0,1\right] $, and for every $\xi
^{1},\xi ^{2}\in 
%TCIMACRO{\U{211d} }%
%BeginExpansion
\mathbb{R}
%EndExpansion
^{nm}$ with $rank\left\{ \xi ^{1}-\xi ^{2}\right\} \leq 1$.

\item[H.2.4] The function $G\left( x,\cdot ,\xi \right) $ is h\"{o}lder
continuous and for almost every $x\in \Omega $ and for every $\xi \in 
%TCIMACRO{\U{211d} }%
%BeginExpansion
\mathbb{R}
%EndExpansion
^{Nn}$ it follows 
\begin{equation}
\left\vert G\left( x,u,\xi \right) -G\left( x,v,\xi \right) \right\vert \leq
c\left( x\right) \left\vert \xi \right\vert ^{\delta }\left\vert
u-v\right\vert ^{\beta }  \label{H.1.4}
\end{equation}%
for a. e. $x\in \Omega $, for every $u,v\in 
%TCIMACRO{\U{211d} }%
%BeginExpansion
\mathbb{R}
%EndExpansion
^{m}$ and for every $\xi \in 
%TCIMACRO{\U{211d} }%
%BeginExpansion
\mathbb{R}
%EndExpansion
^{nm}$ with $0<\delta <q$ and $0<\beta <\min \left\{ \frac{p^{\ast }\left(
q-\delta \right) }{q},\,p-\delta ,1\right\} $, $c$ is a not-negative
function and $c\left( x\right) \in L_{loc}^{\sigma }\left( \Omega \right) $
with $\sigma >\frac{p^{\ast }q}{p^{\ast }\left( q-\delta \right) -\beta q}$
and 
\begin{equation}
\frac{\delta }{p}+\frac{\beta }{p^{\ast }}+\frac{1}{\sigma }<\frac{p}{n}
\label{H.1.4.5}
\end{equation}

\item[H.2.4 (bis)] H.1.4 holds with $c\left( x\right) \in L_{loc}^{\infty
}\left( \Omega \right) $ and 
\begin{equation}
\frac{\delta }{p}+\frac{\beta }{p^{\ast }}<\frac{p}{n}  \label{H.1.4.6}
\end{equation}
\end{description}

\bigskip

IIn particular, the author in [32] proves the following regularity theorem.

\begin{theorem}
If $u\in W_{loc}^{1,p}\left( \Omega ,%
%TCIMACRO{\U{211d} }%
%BeginExpansion
\mathbb{R}
%EndExpansion
^{m}\right) $, with $n\geq 2$, $m\geq 1$ and $1<p<n$, is a local minimizer
of the functional (1.3) and the hypotheses H.2.1, H.2.2, H.2.3 and H.2.4 (or
H.2.4 (bis)) hold then every componets $u^{\alpha }$ of the vectorial
function $u$ are a locally h\"{o}lder continuous functions.
\end{theorem}

\bigskip In [33] the author generalizes these results by considering less
restrictive conditions on the parameters, also considering the polyconvex
case and studying borderline cases. In [31] the author studies the
regularity of the minima of the following class of functionals

\begin{equation}
\mathcal{F}\left( u,\Omega \right) =\int\limits_{\Omega }\sum\limits_{\alpha
=1}^{m}\left\vert \nabla u^{\alpha }\right\vert ^{p}+G\left( x,u,\left\vert
\nabla u^{1}\right\vert ,...,\left\vert \nabla u^{m}\right\vert \right) \,dx
\end{equation}

considering the following hypotheses on density $G$.

\begin{description}
\item[H.3] Let $\Omega $ be a bounded open subset of $%
%TCIMACRO{\U{211d} }%
%BeginExpansion
\mathbb{R}
%EndExpansion
^{n}$ with $n\geq 2\ $and let $G:\Omega \times 
%TCIMACRO{\U{211d} }%
%BeginExpansion
\mathbb{R}
%EndExpansion
^{m}\times 
%TCIMACRO{\U{211d} }%
%BeginExpansion
\mathbb{R}
%EndExpansion
_{0,+}^{m}\rightarrow 
%TCIMACRO{\U{211d} }%
%BeginExpansion
\mathbb{R}
%EndExpansion
$ be a Caratheodory function, where $%
%TCIMACRO{\U{211d} }%
%BeginExpansion
\mathbb{R}
%EndExpansion
_{0,+}=\left[ 0,+\infty \right) $ $\ $and $%
%TCIMACRO{\U{211d} }%
%BeginExpansion
\mathbb{R}
%EndExpansion
_{0,+}^{m}=%
%TCIMACRO{\U{211d} }%
%BeginExpansion
\mathbb{R}
%EndExpansion
_{0,+}\times \cdots \times 
%TCIMACRO{\U{211d} }%
%BeginExpansion
\mathbb{R}
%EndExpansion
_{0,+}$ with $m\geq 1$; we make the following growth conditions on $G$:\
there exists a constant $L>1$ such that%
\begin{equation*}
\sum\limits_{\alpha =1}^{m}\left\vert \xi ^{\alpha }\right\vert
^{q}-\sum\limits_{\alpha =1}^{m}\left\vert s^{\alpha }\right\vert
^{q}-a\left( x\right) \leq G\left( x,s^{1},...,s^{m},\left\vert \xi
^{1}\right\vert ,...,\left\vert \xi ^{m}\right\vert \right) \leq L\left[
\sum\limits_{\alpha =1}^{m}\left\vert \xi ^{\alpha }\right\vert
^{q}+\sum\limits_{\alpha =1}^{m}\left\vert s^{\alpha }\right\vert
^{q}+a\left( x\right) \right]
\end{equation*}%
for $\mathcal{L}^{n}$ a. e. $x\in \Omega $, for every $s^{\alpha }\in 
%TCIMACRO{\U{211d} }%
%BeginExpansion
\mathbb{R}
%EndExpansion
$ and for every $\xi ^{\alpha }\in 
%TCIMACRO{\U{211d} }%
%BeginExpansion
\mathbb{R}
%EndExpansion
$ with $\alpha =1,...,m$ and $m\geq 1$ and with $a\left( x\right) \in
L^{\sigma }\left( \Omega \right) $, $a(x)\geq 0$ for $\mathcal{L}^{n}$ a. e. 
$x\in \Omega $, $\sigma >\frac{n}{p}$, $1\leq q<\frac{p^{2}}{n}$ and $1<p<n$.
\end{description}

\bigskip Assuming that the previous growth hypothesis H.3 holds,the author
in [31] proved the following regularity result.

\begin{theorem}
Let $\Omega $ be a bounded open subset of $%
%TCIMACRO{\U{211d} }%
%BeginExpansion
\mathbb{R}
%EndExpansion
^{n}$ with $n\geq 2$; if $u\in W^{1,p}\left( \Omega ,%
%TCIMACRO{\U{211d} }%
%BeginExpansion
\mathbb{R}
%EndExpansion
^{m}\right) $, with $m\geq 1$, is a local minimum of the functional (1.9)
and $H.3$\ holds then $u^{\alpha }\in C_{loc}^{o,\beta }\left( \Omega
\right) $ for every $\alpha =1,...,m$, with $\beta \in \left( 0,1\right) $.
\end{theorem}

Theorem 1 of this article differs considerably from the previous regularity
results presented in [8-11, 32,33], in fact, unlike the results given in
[8-11, 32,33] no Convexity or Polyconvexity, or Quasi-Convexity, or Rank-one
Convexity hypothesis is made on density $G$, while, unlike the results
contained in [31], density $G$ has no diagonal structure. \ 

In [25] the author with M. Randolfi proved a regularity result for the
minima of vector functionals with anisotropic growths of the following type%
\begin{equation}
\int\limits_{\Omega }\sum\limits_{\alpha =1}^{m}F_{\alpha }\left( x,\nabla
u^{\alpha }\right) +G\left( x,\nabla u\right) \,dx  \label{1.4}
\end{equation}%
with%
\begin{equation}
\sum\limits_{\alpha =1}^{m}\Phi _{i,\alpha }\left( \left\vert \xi
_{i}^{\alpha }\right\vert \right) \leq F_{\alpha }\left( x,\xi ^{\alpha
}\right) \leq L\left[ \bar{B}_{\alpha }^{\beta _{\alpha }}\left( \left\vert
\xi ^{\alpha }\right\vert \right) +a\left( x\right) \right]  \label{1.5}
\end{equation}%
where $\Phi _{i,\alpha }$ are N functions belonging to the class $\triangle
_{2}^{m_{\alpha }}\cap \nabla _{2}^{r_{\alpha }}$, $\bar{B}_{\alpha }$\ is
the Sobolev function associated with $\Phi _{i,\alpha }$'s, $\beta _{\alpha
}\in \left( 0,1\right] $ and $a\in L_{loc}^{\sigma }\left( \Omega \right) $
is a non negative function; oppure with%
\begin{equation}
\sum\limits_{\alpha =1}^{m}\Phi _{i,\alpha }\left( \left\vert \xi
_{i}^{\alpha }\right\vert \right) -a\left( x\right) \leq F_{\alpha }\left(
x,\xi ^{\alpha }\right) \leq L_{1}\left[ \sum\limits_{\alpha =1}^{m}\Phi
_{i,\alpha }\left( \left\vert \xi _{i}^{\alpha }\right\vert \right) +a\left(
x\right) \right]  \label{1.6}
\end{equation}%
where $\Phi _{i,\alpha }$ are N-functions belonging to the class $\triangle
_{2}^{m_{\alpha }}\cap \nabla _{2}^{r_{\alpha }}$ and $a\in L_{loc}^{\sigma
}\left( \Omega \right) $ is a non negative function, moreover, appropriate
hypotheses are made on the density $G$, for more details we refer to [25].
In particular, using the techniques presented in [26-28], the author with M.
Randolfi have shown that the minima of the functional (\ref{1.4}) are
locally bounded functions in the case (\ref{1.5}) and locally Holder
continuous in the case (\ref{1.6}), we refer to [25] for more details.

\bigskip

\section{Preliminary results}

Before giving the proofs of Theorem 1 and Theorem 2, for completeness we
introduce a list of results that we will use during the proof.

\subsection{Lemmata}

\begin{lemma}[Young Inequality]
Let $\varepsilon >0$, $a,b>0$ and $1<p,q<+\infty $ with $\frac{1}{p}+\frac{1%
}{q}=1$\ then it follows%
\begin{equation}
ab\leq \varepsilon \frac{a^{p}}{p}+\frac{b^{q}}{\varepsilon ^{\frac{q}{p}}q}
\end{equation}
\end{lemma}

\begin{lemma}[H\"{o}lder Inequality]
Assume $1\leq p,q\leq +\infty $ with $\frac{1}{p}+\frac{1}{q}=1$\ then if $%
u\in L^{p}\left( \Omega \right) $\ and $v\in L^{p}\left( \Omega \right) $\
it follows%
\begin{equation}
\int\limits_{\Omega }\left\vert uv\right\vert \,dx\leq \left(
\int\limits_{\Omega }\left\vert u\right\vert ^{p}\,dx\right) ^{\frac{1}{p}%
}\left( \int\limits_{\Omega }\left\vert v\right\vert ^{q}\,dx\right) ^{\frac{%
1}{q}}
\end{equation}
\end{lemma}

\begin{lemma}
\label{lem3} Let $Z\left( t\right) $ be a nonnegative and bounded function
on the set $\left[ \varrho ,R\right] $; if for every $\varrho \leq t<s\leq R$
we get%
\begin{equation}
Z\left( t\right) \leq \theta Z\left( s\right) +\frac{A}{\left( s-t\right)
^{\lambda }}+\frac{B}{\left( s-t\right) ^{\mu }}+C
\end{equation}%
where $A,B,C\geq 0$, $\lambda >\mu >0$ and $0\leq \theta <1$ then it follows%
\begin{equation}
Z\left( \varrho \right) \leq C\left( \theta ,\lambda \right) \left( \frac{A}{%
\left( R-\varrho \right) ^{\lambda }}+\frac{B}{\left( R-\varrho \right)
^{\mu }}+C\right)
\end{equation}%
where $C\left( \theta ,\lambda \right) >0$ is a real constant depending only
on $\theta $ and $\lambda $.
\end{lemma}

Refer to [12, 24].

\subsection{Polyconvex, Quasi-Convex and Rank-one Convex functions}

\begin{definition}
\bigskip A function $f:%
%TCIMACRO{\U{211d} }%
%BeginExpansion
\mathbb{R}
%EndExpansion
^{nm}\rightarrow 
%TCIMACRO{\U{211d} }%
%BeginExpansion
\mathbb{R}
%EndExpansion
\cup \left\{ +\infty \right\} $ is said to be rank one convex if%
\begin{equation*}
f\left( \lambda A+\left( 1-\lambda \right) B\right) \leq \lambda f\left(
A\right) +\left( 1-\lambda \right) f\left( B\right)
\end{equation*}%
for every $\lambda \in \left[ 0,1\right] $, $A$,$B\in 
%TCIMACRO{\U{211d} }%
%BeginExpansion
\mathbb{R}
%EndExpansion
^{nm}$ with $rank\left\{ A-B\right\} \leq 1$.
\end{definition}

\begin{definition}
A Borel measurable function and locally integrable function $f:%
%TCIMACRO{\U{211d} }%
%BeginExpansion
\mathbb{R}
%EndExpansion
^{nm}\rightarrow 
%TCIMACRO{\U{211d} }%
%BeginExpansion
\mathbb{R}
%EndExpansion
$ is said to be quasiconvex if%
\begin{equation*}
f\left( A\right) \leq \frac{1}{\left\vert D\right\vert }\int\limits_{D}f%
\left( A+\nabla \varphi \right) \,dx
\end{equation*}%
for every bounded domain $D\subset 
%TCIMACRO{\U{211d} }%
%BeginExpansion
\mathbb{R}
%EndExpansion
^{n}$, for every $A\in 
%TCIMACRO{\U{211d} }%
%BeginExpansion
\mathbb{R}
%EndExpansion
^{nm}$ and for every $\varphi \in W_{0}^{1,\infty }\left( D;%
%TCIMACRO{\U{211d} }%
%BeginExpansion
\mathbb{R}
%EndExpansion
^{nm}\right) $.
\end{definition}

\begin{definition}
A function $f:%
%TCIMACRO{\U{211d} }%
%BeginExpansion
\mathbb{R}
%EndExpansion
^{nm}\rightarrow 
%TCIMACRO{\U{211d} }%
%BeginExpansion
\mathbb{R}
%EndExpansion
\cup \left\{ +\infty \right\} $ is said to be polyconvex if there exists a
function $g:%
%TCIMACRO{\U{211d} }%
%BeginExpansion
\mathbb{R}
%EndExpansion
^{nm}\rightarrow 
%TCIMACRO{\U{211d} }%
%BeginExpansion
\mathbb{R}
%EndExpansion
\cup \left\{ +\infty \right\} $ convex such that%
\begin{equation*}
f\left( A\right) =g\left( T\left( A\right) \right)
\end{equation*}%
where $T:%
%TCIMACRO{\U{211d} }%
%BeginExpansion
\mathbb{R}
%EndExpansion
^{nm}\rightarrow 
%TCIMACRO{\U{211d} }%
%BeginExpansion
\mathbb{R}
%EndExpansion
^{\tau \left( n,m\right) }$ is such that%
\begin{equation*}
T(A)=(A,adj_{2}\left( A\right) ,...,adj_{n\wedge m}\left( A\right) )
\end{equation*}%
where $adj_{s}\left( A\right) $ stands for the matrix of all $s\times s$
minors of tha matrix $A\in 
%TCIMACRO{\U{211d} }%
%BeginExpansion
\mathbb{R}
%EndExpansion
^{nm}$, 2$\leq s\leq n\wedge m=\min \left\{ n,m\right\} $ and%
\begin{equation*}
\tau \left( n,m\right) =\sum\limits_{s=1}^{n\wedge m}\sigma \left( s\right)
\end{equation*}%
where $\sigma \left( s\right) =\frac{n!m!}{\left( s!\right) ^{2}\left(
m-s\right) !\left( n-s\right) !}$.
\end{definition}

In particular we recall the following theorem.

\begin{theorem}

\begin{enumerate}
\item Let $f:%
%TCIMACRO{\U{211d} }%
%BeginExpansion
\mathbb{R}
%EndExpansion
^{nm}\rightarrow 
%TCIMACRO{\U{211d} }%
%BeginExpansion
\mathbb{R}
%EndExpansion
$ then%
\begin{equation*}
f\text{ convex }\Longrightarrow f\text{ polyconvex}\Longrightarrow f\text{
quasiconvex}\Longrightarrow f\text{ rank one convex. }
\end{equation*}

\item If $m=1$ or $n=1$ then all thess notions are equivalent.

\item If $f\in C^{2}\left( 
%TCIMACRO{\U{211d} }%
%BeginExpansion
\mathbb{R}
%EndExpansion
^{nm}\right) $ then rank 0ne convexity is equivalent to Legendre-Hadamard
condition%
\begin{equation*}
\sum\limits_{i,j=1}^{m}\sum\limits_{\alpha ,\beta =1}^{n}\frac{\partial ^{2}f%
}{\partial A_{\alpha }^{i}\partial A_{\beta }^{j}}\left( A\right) \lambda
^{i}\lambda ^{j}\mu _{\alpha }\mu _{\beta }\geq 0
\end{equation*}%
for every $\lambda \in 
%TCIMACRO{\U{211d} }%
%BeginExpansion
\mathbb{R}
%EndExpansion
^{m}$, $\mu \in 
%TCIMACRO{\U{211d} }%
%BeginExpansion
\mathbb{R}
%EndExpansion
^{n}$, $A=\left( A_{\alpha }^{i}\right) _{1\leq i\leq m,1\leq \alpha \leq
n}\in 
%TCIMACRO{\U{211d} }%
%BeginExpansion
\mathbb{R}
%EndExpansion
^{nm}$.

\item If $f:%
%TCIMACRO{\U{211d} }%
%BeginExpansion
\mathbb{R}
%EndExpansion
^{nm}\rightarrow 
%TCIMACRO{\U{211d} }%
%BeginExpansion
\mathbb{R}
%EndExpansion
$ is convex, polyconvex, quasiconvex or rank one convex then $f$ is locally
Lipschitz.
\end{enumerate}
\end{theorem}

Refer to [12, 24].

\subsection{Sobolev Spaces}

\begin{theorem}[Sobolev Inequality]
Let $\Omega $ be a open subset of $%
%TCIMACRO{\U{211d} }%
%BeginExpansion
\mathbb{R}
%EndExpansion
^{N}$ if $u\in W_{0}^{1,p}\left( \Omega \right) $ with $1\leq p<N$ there
exists a real positive constant $C_{SN}$, depending only on $p$ and $N$,
such that%
\begin{equation}
\left\Vert u\right\Vert _{L^{p^{\ast }}\left( \Omega \right) }\leq
C_{SN}\left\Vert \nabla u\right\Vert _{L^{p}\left( \Omega \right) }
\end{equation}%
where $p^{\ast }=\frac{Np}{N-p}$.
\end{theorem}

\begin{theorem}
(Rellich-Sobolev Immersion Theorem) Let $\Omega $ be a open bounded subset
of $%
%TCIMACRO{\U{211d} }%
%BeginExpansion
\mathbb{R}
%EndExpansion
^{N}$ with lipschitz boundary then if $u\in W^{1,p}\left( \Omega \right) $
with $1\leq p<N$ there exists a real positive constant $C_{IS}$, depending
only on $p$ and $N$, such that%
\begin{equation}
\left\Vert u\right\Vert _{L^{p^{\ast }}\left( \Omega \right) }\leq
C_{IS}\left\Vert u\right\Vert _{W^{1,p}\left( \Omega \right) }
\end{equation}%
where $p^{\ast }=\frac{Np}{N-p}$.
\end{theorem}

Refer to [3, 12, 24, 40, 41].

For completeness we remember that if $\Omega $ is a open subset of $%
%TCIMACRO{\U{211d} }%
%BeginExpansion
\mathbb{R}
%EndExpansion
^{N}$ and $u$\ is a Lebesgue measurable function then $L^{p}\left( \Omega
\right) $ is the set of the class of the Lebesgue measurable function such
that $\int\limits_{\Omega }\left\vert u\right\vert ^{p}\,dx<+\infty $ and $%
W^{1,p}\left( \Omega \right) $\ is the set of the function $u\in L^{p}\left(
\Omega \right) $ such that its waek derivate $\partial _{i}u\in L^{p}\left(
\Omega \right) $. The spaces $L^{p}\left( \Omega \right) $ and $%
W^{1,p}\left( \Omega \right) $ are Banach spaces with the respective norms 
\begin{equation}
\left\Vert u\right\Vert _{L^{p}\left( \Omega \right) }=\left(
\int\limits_{\Omega }\left\vert u\right\vert ^{p}\,dx\right) ^{\frac{1}{p}}
\end{equation}%
and%
\begin{equation}
\left\Vert u\right\Vert _{W^{1,p}\left( \Omega \right) }=\left\Vert
u\right\Vert _{L^{p}\left( \Omega \right) }+\sum\limits_{i=1}^{N}\left\Vert
\partial _{i}u\right\Vert _{L^{p}\left( \Omega \right) }
\end{equation}%
We say that the function $u:\Omega \subset 
%TCIMACRO{\U{211d} }%
%BeginExpansion
\mathbb{R}
%EndExpansion
^{N}\rightarrow 
%TCIMACRO{\U{211d} }%
%BeginExpansion
\mathbb{R}
%EndExpansion
^{n}$ belong in $W^{1,p}\left( \Omega ,%
%TCIMACRO{\U{211d} }%
%BeginExpansion
\mathbb{R}
%EndExpansion
^{n}\right) $ if $u^{\alpha }\in W^{1,p}\left( \Omega \right) $ for every $%
\alpha =1,...,n$, where $u^{\alpha }$ is the $\alpha $ component of the
vector-valued function $u$; we end by remembering that $W^{1,p}\left( \Omega
,%
%TCIMACRO{\U{211d} }%
%BeginExpansion
\mathbb{R}
%EndExpansion
^{n}\right) $ is a Banach space with the norm%
\begin{equation}
\left\Vert u\right\Vert _{W^{1,p}\left( \Omega ,%
%TCIMACRO{\U{211d} }%
%BeginExpansion
\mathbb{R}
%EndExpansion
^{n}\right) }=\sum\limits_{\alpha =1}^{n}\left\Vert u^{\alpha }\right\Vert
_{W^{1,p}\left( \Omega \right) }
\end{equation}

\begin{definition}
Let $\Omega \subset 
%TCIMACRO{\U{211d} }%
%BeginExpansion
\mathbb{R}
%EndExpansion
^{N}$ be a bounded open set and $v:\Omega \rightarrow 
%TCIMACRO{\U{211d} }%
%BeginExpansion
\mathbb{R}
%EndExpansion
$, we say that $v\in W_{loc}^{1,p}\left( \Omega \right) $ belong to the De
Giorgi class $DG^{+}\left( \Omega ,p,\lambda ,\lambda _{\ast },\chi
,\varepsilon ,R_{0},k_{0}\right) $ with $p>1$, $\lambda >0$, $\lambda _{\ast
}>0$, $\chi >0$, $\varepsilon >0$, $R_{0}>0$ and $k_{0}\geq 0$ if%
\begin{equation}
\int\limits_{A_{k,\varrho }}\left\vert \nabla v\right\vert ^{p}\,dx\leq 
\frac{\lambda }{\left( R-\varrho \right) ^{p}}\int\limits_{A_{k,R}}\left(
v-k\right) ^{p}\,dx+\lambda _{\ast }\left( \chi ^{p}+k^{p}R^{-N\varepsilon
}\right) \left\vert A_{k,R}\right\vert ^{1-\frac{p}{N}+\varepsilon }
\end{equation}%
for all $k\geq k_{0}\geq 0$ and for all pair of balls $B_{\varrho }\left(
x_{0}\right) \subset B_{R}\left( x_{0}\right) \subset \subset \Omega $ with $%
0<\varrho <R<R_{0}$ and $A_{k,s}=B_{s}\left( x_{0}\right) \cap \left\{
v>k\right\} $ with $s>0$.
\end{definition}

\begin{definition}
Let $\Omega \subset 
%TCIMACRO{\U{211d} }%
%BeginExpansion
\mathbb{R}
%EndExpansion
^{N}$ be a bounded open set and $v:\Omega \rightarrow 
%TCIMACRO{\U{211d} }%
%BeginExpansion
\mathbb{R}
%EndExpansion
$, we say that $v\in W_{loc}^{1,p}\left( \Omega \right) $ belong to the De
Giorgi class $DG^{-}\left( \Omega ,p,\lambda ,\lambda _{\ast },\chi
,\varepsilon ,R_{0},k_{0}\right) $ with $p>1$, $\lambda >0$, $\lambda _{\ast
}>0$, $\chi >0$ and $k_{0}\geq 0$ if%
\begin{equation}
\int\limits_{B_{k,\varrho }}\left\vert \nabla v\right\vert ^{p}\,dx\leq 
\frac{\lambda }{\left( R-\varrho \right) ^{p}}\int\limits_{B_{k,R}}\left(
k-v\right) ^{p}\,dx+\lambda _{\ast }\left( \chi ^{p}+\left\vert k\right\vert
^{p}R^{-N\varepsilon }\right) \left\vert B_{k,R}\right\vert ^{1-\frac{p}{N}%
+\varepsilon }
\end{equation}%
for all $k\leq -k_{0}\leq 0$ and for all pair of balls $B_{\varrho }\left(
x_{0}\right) \subset B_{R}\left( x_{0}\right) \subset \subset \Omega $ with $%
0<\varrho <R<R_{0}$ and $B_{k,s}=B_{s}\left( x_{0}\right) \cap \left\{
v<k\right\} $ with $s>0$.
\end{definition}

\begin{definition}
We set $DG\left( \Omega ,p,\lambda ,\lambda _{\ast },\chi ,\varepsilon
,R_{0},k_{0}\right) =DG^{+}\left( \Omega ,p,\lambda ,\lambda _{\ast },\chi
,\varepsilon ,R_{0},k_{0}\right) \cap DG^{-}\left( \Omega ,p,\lambda
,\lambda _{\ast },\chi ,\varepsilon ,R_{0},k_{0}\right) $.
\end{definition}

\begin{theorem}
Let $v\in DG\left( \Omega ,p,\lambda ,\lambda _{\ast },\chi ,\varepsilon
,R_{0},k_{0}\right) $ and $\tau \in (0,1)$, then there exists a constant $%
C>1 $ depending only upon the data and not-dependent on $v$ and $x_{0}\in
\Omega $ such that for every pair of balls $B_{\tau \varrho }\left(
x_{0}\right) \subset B_{\varrho }\left( x_{0}\right) \subset \subset \Omega $
with $0<\varrho <R_{0}$ 
\begin{equation}
\left\Vert v\right\Vert _{L^{\infty }\left( B_{\tau \varrho }\left(
x_{0}\right) \right) }\leq \max \left\{ \lambda _{\ast }\varrho ^{\frac{%
N\varepsilon }{p}};\frac{C}{\left( 1-\tau \right) ^{\frac{N}{p}}}\left[ 
\frac{1}{\left\vert B_{\varrho }\left( x_{0}\right) \right\vert }%
\int\limits_{B_{\varrho }\left( x_{0}\right) }\left\vert v\right\vert
^{p}\,dx\right] ^{\frac{1}{p}}\right\}
\end{equation}%
moreover, there exists $\tilde{\alpha}\in (0,1)$ depending only upon the
data and not-dependent on $v$ and $x_{0}\in \Omega $ such that 
\begin{equation}
osc(v,B_{\varrho }\left( x_{0}\right) )\leq C\max \left\{ \lambda _{\ast
}\varrho ^{\frac{N\varepsilon }{p}};\left( \frac{\varrho }{R}\right) ^{%
\tilde{\alpha}}osc(v,B_{R}\left( x_{0}\right) )\right\}
\end{equation}%
where $osc(v,B_{s}\left( x_{0}\right) )=ess\,\sup_{B_{s}\left( x_{0}\right)
}\left( v\right) -ess$ $inf_{B_{s}\left( x_{0}\right) }\left( v\right) $.
Therefore $v\in C_{loc}^{0,\tilde{\alpha}_{0}}\left( \Omega \right) $ with $%
\tilde{\alpha}_{0}=\max \left\{ \tilde{\alpha};N\chi \right\} $.
\end{theorem}

For more details on De Giorgi's classes and for the proof of the Theorem 8
refer to [22, 24] (see olso [13, 38, 39] for the De Giorgi--Moser-Nash
Theorem).

\section{The proof of Theorem 2}

Let us consider $y\in \Omega $ then we fix $R_{0}=\frac{1}{4}\min \left\{ 
\frac{1}{\sqrt[N]{\varpi _{N}}},dist\left( \partial \Omega ,y\right)
\right\} $, where\ $\varpi _{N}=\left\vert B_{1}\left( 0\right) \right\vert $%
, and we define $\Sigma =\left\{ x\in \Omega :\left\vert x-y\right\vert \leq
R_{0}\right\} $. We fix $x_{0}\in \Sigma $, $R_{1}=\frac{1}{4}dist\left(
\partial \Sigma ,x_{0}\right) $, $0<\varrho \leq t<s\leq R<R_{0}$, $%
B_{z}\left( x_{0}\right) =\left\{ x:\left\vert x-x_{0}\right\vert <z\right\} 
$\ and we choose $\eta \in C_{c}^{\infty }\left( B_{s}\left( x_{0}\right)
\right) $\ such that $\eta =1$\ on $B_{t}\left( x_{0}\right) $, $0\leq \eta
\leq 1$\ on $B_{s}\left( x_{0}\right) $\ and $\left\vert \nabla \eta
\right\vert \leq \frac{2}{s-t}$\ on $B_{s}\left( x_{0}\right) $. Let us
define%
\begin{equation*}
\varphi =-\eta ^{p}w
\end{equation*}%
where $w\in W^{1,p}\left( \Sigma ,%
%TCIMACRO{\U{211d} }%
%BeginExpansion
\mathbb{R}
%EndExpansion
^{n}\right) $ with%
\begin{equation*}
w^{1}=\max \left( u^{1}-k,0\right) ,w^{\alpha }=0,\alpha =2,...,n
\end{equation*}%
Let us observe that $\varphi =0$ $\mathcal{L}^{N}$-a.e. in $\Omega
\backslash \left( \left\{ \eta >0\right\} \cap \left\{ u^{1}>k\right\}
\right) $ thus%
\begin{equation}
\nabla u+\nabla \varphi =\nabla u
\end{equation}%
$\mathcal{L}^{N}$-a.e. in $\Omega \backslash \left( \left\{ \eta >0\right\}
\cap \left\{ u^{1}>k\right\} \right) $. Since $u$ is a local minimizer of
the functional (1.1) then we get%
\begin{equation}
J\left( u,\Sigma \right) \leq J\left( u+\varphi ,\Sigma \right)
\end{equation}%
it is%
\begin{equation}
\begin{tabular}{l}
$\int\limits_{\Sigma }\sum\limits_{\alpha =1}^{n}\left\vert \nabla u^{\alpha
}\right\vert ^{p}+G\left( x,u,\nabla u\right) \,dx$ \\ 
$\leq \int\limits_{\Sigma }\sum\limits_{\alpha =1}^{n}\left\vert \nabla
u^{\alpha }+\nabla \varphi ^{\alpha }\right\vert ^{p}+G\left( x,u+\varphi
,\nabla u+\nabla \varphi \right) \,dx$%
\end{tabular}%
\end{equation}%
and%
\begin{equation}
\begin{tabular}{l}
$\int\limits_{\Sigma }\sum\limits_{\alpha =2}^{n}\left\vert \nabla u^{\alpha
}\right\vert ^{p}\,dx+\int\limits_{\Sigma }\left\vert \nabla
u^{1}\right\vert ^{p}+G\left( x,u,\nabla u\right) \,dx$ \\ 
$\leq \int\limits_{\Sigma }\sum\limits_{\alpha =2}^{n}\left\vert \nabla
u^{\alpha }\right\vert ^{p}\,dx+\int\limits_{\Sigma }\left\vert \nabla
u^{1}+\nabla \varphi ^{1}\right\vert ^{p}+G\left( x,u+\varphi ,\nabla
u+\nabla \varphi \right) \,dx$%
\end{tabular}%
\end{equation}%
From (2.20) we deduce%
\begin{equation}
\begin{tabular}{l}
$\int\limits_{\Sigma }\left\vert \nabla u^{1}\right\vert ^{p}+G\left(
x,u,\nabla u\right) \,dx$ \\ 
$\leq \int\limits_{\Sigma }\left\vert \nabla u^{1}+\nabla \varphi
^{1}\right\vert ^{p}+G\left( x,u+\varphi ,\nabla u+\nabla \varphi \right)
\,dx$ \\ 
$=\int\limits_{B_{r}\left( x_{0}\right) }\left\vert \nabla u^{1}+\nabla
\varphi ^{1}\right\vert ^{p}\,dx+\int\limits_{\Sigma -B_{r}\left(
x_{0}\right) \backslash }\left\vert \nabla u^{1}\right\vert ^{p}\,dx$ \\ 
$+\int\limits_{B_{s}\left( x_{0}\right) }G\left( x,u+\varphi ,\nabla
u+\nabla \varphi \right) \,dx+\int\limits_{\Sigma -B_{s}\left( x_{0}\right)
\backslash }G\left( x,u,\nabla u\right) \,dx$%
\end{tabular}%
\end{equation}%
and%
\begin{equation}
\begin{tabular}{l}
$\int\limits_{B_{s}\left( x_{0}\right) }\left\vert \nabla u^{1}\right\vert
^{p}+G\left( x,u,\nabla u\right) \,dx$ \\ 
$\leq \int\limits_{B_{r}\left( x_{0}\right) }\left\vert \nabla u^{1}+\nabla
\varphi ^{1}\right\vert ^{p}\,dx+$ \\ 
$+\int\limits_{B_{s}\left( x_{0}\right) }G\left( x,u+\varphi ,\nabla
u+\nabla \varphi \right) \,dx$%
\end{tabular}%
\end{equation}%
Let us define $E_{k,s}^{1}=\left\{ \eta >0\right\} \cap \left\{
u^{1}>k\right\} \cap B_{s}\left( x_{0}\right) \subset B_{s}\left(
x_{0}\right) $\ then%
\begin{equation}
\begin{tabular}{l}
$\int\limits_{E_{k,s}^{1}}\left\vert \nabla u^{1}\right\vert ^{p}+G\left(
x,u,\nabla u\right) \,dx+\int\limits_{B_{s}\left( x_{0}\right)
-E_{k,s}^{1}}\left\vert \nabla u^{1}\right\vert ^{p}+G\left( x,u,\nabla
u\right) \,dx$ \\ 
$\leq \int\limits_{E_{k,s}^{1}}\left\vert \nabla u^{1}+\nabla \varphi
^{1}\right\vert ^{p}\,dx+\int\limits_{B_{s}\left( x_{0}\right)
-E_{k,s}^{1}}\left\vert \nabla u^{1}\right\vert ^{p}\,dx$ \\ 
$+\int\limits_{E_{k,s}^{1}}G\left( x,u+\varphi ,\nabla u+\nabla \varphi
\right) \,dx$ \\ 
$+\int\limits_{B_{s}\left( x_{0}\right) -E_{k,s}^{1}}G\left( x,u,\nabla
u\right) \,dx$%
\end{tabular}%
\end{equation}%
and%
\begin{equation}
\begin{tabular}{l}
$\int\limits_{E_{k,s}^{1}}\left\vert \nabla u^{1}\right\vert ^{p}\,dx$ \\ 
$\leq \int\limits_{E_{k,s}^{1}}\left\vert \nabla u^{1}+\nabla \varphi
^{1}\right\vert ^{p}\,dx+$ \\ 
$+\int\limits_{E_{k,s}^{1}}G\left( x,u+\varphi ,\nabla u+\nabla \varphi
\right) -G\left( x,u,\nabla u\right) \,dx$%
\end{tabular}
\label{45}
\end{equation}%
Since 
\begin{equation}
\begin{tabular}{l}
$\int\limits_{E_{k,s}^{1}}G\left( x,u+\varphi ,\nabla u+\nabla \varphi
\right) -G\left( x,u,\nabla u\right) \,dx$ \\ 
$\leq \int\limits_{E_{k,s}^{1}}\left\vert G\left( x,u+\varphi ,\nabla
u+\nabla \varphi \right) -G\left( x,u,\nabla u\right) \right\vert \,dx$%
\end{tabular}
\label{46}
\end{equation}%
using (\ref{45}), (\ref{46}) and H.1.1 we get%
\begin{equation}
\begin{tabular}{l}
$\int\limits_{E_{k,s}^{1}}\left\vert \nabla u^{1}\right\vert ^{p}\,dx$ \\ 
$\leq \int\limits_{E_{k,s}^{1}}\left\vert \nabla u^{1}+\nabla \varphi
^{1}\right\vert ^{p}\,dx+$ \\ 
$+\int\limits_{E_{k,s}^{1}}a\left( x\right) \left\vert \varphi \right\vert
^{\alpha }\left( \left\vert \nabla u+\nabla \varphi \right\vert +\left\vert
\nabla u\right\vert +1\right) ^{q_{1}}+b\left( x\right) \left( \left\vert
u+\varphi \right\vert +\left\vert u\right\vert +1\right) ^{q_{2}}\left\vert
\nabla \varphi \right\vert ^{\beta }\,dx$%
\end{tabular}%
\end{equation}%
Since%
\begin{equation}
\begin{tabular}{l}
$\int\limits_{E_{k,s}^{1}}\left\vert \nabla u^{1}+\nabla \varphi
^{1}\right\vert ^{p}\,dx$ \\ 
$\leq 2^{p-1}\int\limits_{E_{k,s}^{1}}\left( 1-\eta ^{p}\right) \left\vert
\nabla u^{1}\right\vert ^{p}\,dx$ \\ 
$+2^{2p-1}p^{p}\int\limits_{E_{k,s}^{1}-E_{k,t}^{1}}\frac{\left(
u^{1}-k\right) ^{p}}{\left( s-t\right) ^{p}}\,dx$%
\end{tabular}
\label{54}
\end{equation}%
it follows%
\begin{equation}
\begin{tabular}{l}
$\int\limits_{E_{k,s}^{1}}\left\vert \nabla u^{1}\right\vert ^{p}\,dx$ \\ 
$\leq 2^{p-1}\int\limits_{E_{k,s}^{1}}\left( 1-\eta ^{p}\right) \left\vert
\nabla u^{1}\right\vert
^{p}\,dx++2^{2p-1}p^{p}\int\limits_{E_{k,s}^{1}-E_{k,t}^{1}}\frac{\left(
u^{1}-k\right) ^{p}}{\left( s-t\right) ^{p}}\,dx$ \\ 
$+\int\limits_{E_{k,s}^{1}}a\left( x\right) \left\vert \varphi \right\vert
^{\alpha }\left( \left\vert \nabla u+\nabla \varphi \right\vert +\left\vert
\nabla u\right\vert +1\right) ^{q_{1}}+b\left( x\right) \left( \left\vert
u+\varphi \right\vert +\left\vert u\right\vert +1\right) ^{q_{2}}\left\vert
\nabla \varphi \right\vert ^{\beta }\,dx$%
\end{tabular}%
\end{equation}

\bigskip Now let's estimate the following term%
\begin{equation*}
\int\limits_{E_{k,s}^{1}}a\left( x\right) \left\vert \varphi \right\vert
^{\alpha }\left( \left\vert \nabla u+\nabla \varphi \right\vert +\left\vert
\nabla u\right\vert +1\right) ^{q_{1}}+b\left( x\right) \left( \left\vert
u+\varphi \right\vert +\left\vert u\right\vert +1\right) ^{q_{2}}\left\vert
\nabla \varphi \right\vert ^{\beta }\,dx
\end{equation*}%
since $\alpha ,\beta <1$, using H\"{o}lder's inequality, we obtain%
\begin{equation}
\begin{tabular}{l}
$\int\limits_{E_{k,s}^{1}}a\left( x\right) \left\vert \varphi \right\vert
^{\alpha }\left( \left\vert \nabla u+\nabla \varphi \right\vert +\left\vert
\nabla u\right\vert +1\right) ^{q_{1}}+b\left( x\right) \left( \left\vert
u+\varphi \right\vert +\left\vert u\right\vert +1\right) ^{q_{2}}\left\vert
\nabla \varphi \right\vert ^{\beta }\,dx$ \\ 
$\leq \int\limits_{E_{k,s}^{1}}a\left( x\right) \left\vert \varphi
\right\vert ^{\alpha }\left( \left\vert \nabla \varphi \right\vert
+2\left\vert \nabla u\right\vert +1\right) ^{q_{1}}+b\left( x\right) \left(
\left\vert \varphi \right\vert +2\left\vert u\right\vert +1\right)
^{q_{2}}\left\vert \nabla \varphi \right\vert ^{\beta }\,dx$ \\ 
$\leq \left\Vert \varphi \right\Vert _{L^{p^{\ast }}\left( B_{s}\left(
x_{0}\right) \right) }^{\alpha }\left( \int\limits_{E_{k,s}^{1}}\left(
a\left( x\right) \right) ^{\frac{p^{\ast }}{p^{\ast }-\alpha }}\left(
\left\vert \nabla \varphi \right\vert +2\left\vert \nabla u\right\vert
+1\right) ^{\frac{q_{1}p^{\ast }}{p^{\ast }-\alpha }}\,dx\right) ^{\frac{%
p^{\ast }-\alpha }{p^{\ast }}}$ \\ 
$+\left( \int\limits_{E_{k,s}^{1}}\left\vert \nabla \varphi \right\vert
^{p}\,dx\right) ^{\frac{\beta }{p}}\left( \int\limits_{E_{k,s}^{1}}\left(
b\left( x\right) \right) ^{\frac{p}{p-\beta }}\left( \left\vert \varphi
\right\vert +2\left\vert u\right\vert +1\right) ^{\frac{pq_{2}}{p-\beta }%
}\,dx\right) ^{\frac{p-\beta }{p}}$%
\end{tabular}%
\end{equation}%
moreover, since $\alpha <\min \left\{ 1,p-q_{1}\right\} <\frac{p^{\ast
}\left( p-q_{1}\right) }{p}$ and $\beta <\min \left\{ 1,p-q_{2}\right\} <%
\frac{p\left( p^{\ast }-q_{2}\right) }{p^{\ast }}$, then using H\"{o}lder's
inequality and Sobolev's inequality it follows that%
\begin{equation*}
\begin{tabular}{l}
$\int\limits_{E_{k,s}^{1}}a\left( x\right) \left\vert \varphi \right\vert
^{\alpha }\left( \left\vert \nabla u+\nabla \varphi \right\vert +\left\vert
\nabla u\right\vert +1\right) ^{q_{1}}+b\left( x\right) \left( \left\vert
u+\varphi \right\vert +\left\vert u\right\vert +1\right) ^{q_{2}}\left\vert
\nabla \varphi \right\vert ^{\beta }\,dx$ \\ 
$\leq C_{SI}^{\alpha }\left( \int\limits_{E_{k,s}^{1}}\left\vert \nabla
\varphi \right\vert ^{p}\,dx\right) ^{\frac{\alpha }{p}}$ \\ 
$\cdot \left( \left( \int\limits_{E_{k,s}^{1}}\left( a\left( x\right)
\right) ^{\frac{pp^{\ast }}{p\left( p^{\ast }-\alpha \right) -q_{1}p^{\ast }}%
}\,dx\right) ^{\frac{p\left( p^{\ast }-\alpha \right) -q_{1}p^{\ast }}{%
p\left( p^{\ast }-\alpha \right) }}\left( \int\limits_{E_{k,s}^{1}}\left(
\left\vert \nabla \varphi \right\vert +2\left\vert \nabla u\right\vert
+1\right) ^{p}\,dx\right) ^{\frac{q_{1}p^{\ast }}{\left( p^{\ast }-\alpha
\right) p}}\right) ^{\frac{p^{\ast }-\alpha }{p^{\ast }}}$ \\ 
$+\left( \int\limits_{E_{k,s}^{1}}\left\vert \nabla \varphi \right\vert
^{p}\,dx\right) ^{\frac{\beta }{p}}\left( \left(
\int\limits_{E_{k,s}^{1}}\left( b\left( x\right) \right) ^{\frac{pp^{\ast }}{%
\left( p-\beta \right) p^{\ast }-pq_{2}}}\,dx\right) ^{\frac{\left( p-\beta
\right) p^{\ast }-pq_{2}}{\left( p-\beta \right) p^{\ast }}}\left(
\int\limits_{E_{k,s}^{1}}\left( \left\vert \varphi \right\vert +2\left\vert
u\right\vert +1\right) ^{p^{\ast }}\,dx\right) ^{\frac{pq_{2}}{p^{\ast
}(p-\beta )}}\right) ^{\frac{p-\beta }{p}}$%
\end{tabular}%
\end{equation*}%
Now, since $\sigma _{1}>\frac{pp^{\ast }}{\left( p^{\ast }-\alpha \right)
p-q_{1}p^{\ast }}$ and $\sigma _{2}>\frac{pp^{\ast }}{\left( p-\beta \right)
p^{\ast }-q_{2}p}$ we get%
\begin{equation}
\begin{tabular}{l}
$\int\limits_{E_{k,s}^{1}}a\left( x\right) \left\vert \varphi \right\vert
^{\alpha }\left( \left\vert \nabla u+\nabla \varphi \right\vert +\left\vert
\nabla u\right\vert +1\right) ^{q_{1}}+b\left( x\right) \left( \left\vert
u+\varphi \right\vert +\left\vert u\right\vert +1\right) ^{q_{2}}\left\vert
\nabla \varphi \right\vert ^{\beta }\,dx$ \\ 
$\leq C_{SI}^{\alpha }\left( \int\limits_{E_{k,s}^{1}}\left\vert \nabla
\varphi \right\vert ^{p}\,dx\right) ^{\frac{\alpha }{p}}$ \\ 
$\cdot \left( \left( \left\vert E_{k,s}^{1}\right\vert ^{1-\frac{pp^{\ast }}{%
\left[ p\left( p^{\ast }-\alpha \right) -q_{1}p^{\ast }\right] \sigma _{1}}%
}\left\Vert a\right\Vert _{L^{\sigma _{1}}\left( E_{k,s}^{1}\right) }^{\frac{%
pp^{\ast }}{\left[ p\left( p^{\ast }-\alpha \right) -q_{1}p^{\ast }\right] }%
}\right) ^{\frac{p\left( p^{\ast }-\alpha \right) -q_{1}p^{\ast }}{p\left(
p^{\ast }-\alpha \right) }}\left( \int\limits_{E_{k,s}^{1}}\left( \left\vert
\nabla \varphi \right\vert +2\left\vert \nabla u\right\vert +1\right)
^{p}\,dx\right) ^{\frac{q_{1}p^{\ast }}{\left( p^{\ast }-\alpha \right) p}%
}\right) ^{\frac{p^{\ast }-\alpha }{p^{\ast }}}$ \\ 
$+\left( \int\limits_{E_{k,s}^{1}}\left\vert \nabla \varphi \right\vert
^{p}\,dx\right) ^{\frac{\beta }{p}}$ \\ 
$\cdot \left( \left( \left\vert E_{k,s}^{1}\right\vert ^{1-\frac{pp^{\ast }}{%
\left[ \left( p-\beta \right) p^{\ast }-pq_{2}\right] \sigma _{2}}%
}\left\Vert b\right\Vert _{L^{\sigma _{2}}\left( E_{k,s}^{1}\right) }^{\frac{%
pp^{\ast }}{\left( p-\beta \right) p^{\ast }-pq_{2}}}\right) ^{\frac{\left(
p-\beta \right) p^{\ast }-pq_{2}}{\left( p-\beta \right) p^{\ast }}}\left(
\int\limits_{E_{k,s}^{1}}\left( \left\vert \varphi \right\vert +2\left\vert
u\right\vert +1\right) ^{p^{\ast }}\,dx\right) ^{\frac{pq_{2}}{p^{\ast
}(p-\beta )}}\right) ^{\frac{p-\beta }{p}}$ \\ 
$\leq \left( \int\limits_{E_{k,s}^{1}}\left\vert \nabla \varphi \right\vert
^{p}\,dx\right) ^{\frac{\alpha }{p}}\left( \int\limits_{E_{k,s}^{1}}\left(
\left\vert \nabla \varphi \right\vert +2\left\vert \nabla u\right\vert
+1\right) ^{p}\,dx\right) ^{\frac{q_{1}}{p}}\left( C_{SI}^{\alpha
}\left\vert E_{k,s}^{1}\right\vert ^{\left( \frac{p\left( p^{\ast }-\alpha
\right) -q_{1}p^{\ast }}{pp^{\ast }}-\frac{1}{\sigma _{1}}\right)
}\left\Vert a\right\Vert _{L^{\sigma _{1}}\left( E_{k,s}^{1}\right) }\right) 
$ \\ 
$+\left( \int\limits_{E_{k,s}^{1}}\left\vert \nabla \varphi \right\vert
^{p}\,dx\right) ^{\frac{\beta }{p}}\left( \int\limits_{E_{k,s}^{1}}\left(
\left\vert \varphi \right\vert +2\left\vert u\right\vert +1\right) ^{p^{\ast
}}\,dx\right) ^{\frac{q_{2}}{p^{\ast }}}\left( \left\vert
E_{k,s}^{1}\right\vert ^{\frac{\left( p-\beta \right) p^{\ast }-pq_{2}}{%
pp^{\ast }}-\frac{1}{\sigma _{2}}}\left\Vert b\right\Vert _{L^{\sigma
_{2}}\left( E_{k,s}^{1}\right) }\right) $%
\end{tabular}%
\end{equation}%
and%
\begin{equation*}
\begin{tabular}{l}
$\int\limits_{E_{k,s}^{1}}a\left( x\right) \left\vert \varphi \right\vert
^{\alpha }\left( \left\vert \nabla u+\nabla \varphi \right\vert +\left\vert
\nabla u\right\vert +1\right) ^{q_{1}}+b\left( x\right) \left( \left\vert
u+\varphi \right\vert +\left\vert u\right\vert +1\right) ^{q_{2}}\left\vert
\nabla \varphi \right\vert ^{\beta }\,dx$ \\ 
$\leq \left( \int\limits_{E_{k,s}^{1}}\left\vert \nabla \varphi \right\vert
^{p}\,dx\right) ^{\frac{\alpha }{p}}\left( \int\limits_{E_{k,s}^{1}}\left(
\left\vert \nabla \varphi \right\vert +2\left\vert \nabla u\right\vert
+1\right) ^{p}\,dx\right) ^{\frac{q_{1}}{p}}\left( C_{SI}^{\alpha
}\left\vert E_{k,s}^{1}\right\vert ^{\Theta _{1}}\left\Vert a\right\Vert
_{L^{\sigma _{1}}\left( E_{k,s}^{1}\right) }\right) $ \\ 
$+\left( \int\limits_{E_{k,s}^{1}}\left\vert \nabla \varphi \right\vert
^{p}\,dx\right) ^{\frac{\beta }{p}}\left( \int\limits_{E_{k,s}^{1}}\left(
\left\vert \varphi \right\vert +2\left\vert u\right\vert +1\right) ^{p^{\ast
}}\,dx\right) ^{\frac{q_{2}}{p^{\ast }}}\left( \left\vert
E_{k,s}^{1}\right\vert ^{\Theta _{2}}\left\Vert b\right\Vert _{L^{\sigma
_{2}}\left( E_{k,s}^{1}\right) }\right) $%
\end{tabular}%
\end{equation*}%
where $\Theta _{1}=1-\left( \frac{q_{1}}{p}+\frac{\alpha }{p^{\ast }}+\frac{1%
}{\sigma _{1}}\right) $ and $\Theta _{2}=1-\left( \frac{\beta }{p}+\frac{%
q_{2}}{p^{\ast }}+\frac{1}{\sigma _{2}}\right) $. Since%
\begin{equation}
\begin{tabular}{l}
$\left( \int\limits_{E_{k,s}^{1}}\left\vert \nabla \varphi \right\vert
^{p}\,dx\right) ^{\frac{\alpha }{p}}\left( \int\limits_{E_{k,s}^{1}}\left(
\left\vert \nabla \varphi \right\vert +2\left\vert \nabla u\right\vert
+1\right) ^{p}\,dx\right) ^{\frac{q_{1}}{p}}\left( C_{SI}^{\alpha
}\left\vert E_{k,s}^{1}\right\vert ^{\Theta _{1}}\left\Vert a\right\Vert
_{L^{\sigma _{1}}\left( E_{k,s}^{1}\right) }\right) $ \\ 
$\leq 2^{\left( \frac{p-1}{p}\right) q_{1}}\left( C_{SI}^{\alpha }\left\vert
E_{k,s}^{1}\right\vert ^{\Theta _{1}}\left\Vert a\right\Vert _{L^{\sigma
_{1}}\left( E_{k,s}^{1}\right) }\right) \left(
\int\limits_{E_{k,s}^{1}}\left\vert \nabla \varphi \right\vert
^{p}\,dx\right) ^{\frac{\alpha }{p}}$ \\ 
$\cdot \left[ \left( \int\limits_{E_{k,s}^{1}}\left\vert \nabla \varphi
\right\vert ^{p}\,dx\right) ^{\frac{q_{1}}{p}}+2^{\left( \frac{2p-1}{p}%
\right) q_{1}}\left( \int\limits_{E_{k,s}^{1}}\left\vert \nabla u\right\vert
^{p}\,dx\right) ^{\frac{q_{1}}{p}}+2^{\left( \frac{2p-1}{p}\right)
q_{1}}\left\vert E_{k,s}^{1}\right\vert ^{\frac{q_{1}}{p}}\right] $ \\ 
$\leq 2^{\left( \frac{p-1}{p}\right) q_{1}}\left( C_{SI}^{\alpha }\left\vert
E_{k,s}^{1}\right\vert ^{\Theta _{1}}\left\Vert a\right\Vert _{L^{\sigma
_{1}}\left( E_{k,s}^{1}\right) }\right) $ \\ 
$\cdot \left[ \left( \int\limits_{E_{k,s}^{1}}\left\vert \nabla \varphi
\right\vert ^{p}\,dx\right) ^{\frac{\alpha +q_{1}}{p}}+2^{\left( \frac{2p-1}{%
p}\right) q_{1}}\left( \int\limits_{E_{k,s}^{1}}\left\vert \nabla
u\right\vert ^{p}\,dx\right) ^{\frac{q_{1}}{p}}\left(
\int\limits_{E_{k,s}^{1}}\left\vert \nabla \varphi \right\vert
^{p}\,dx\right) ^{\frac{\alpha }{p}}\right. $ \\ 
$\left. +2^{\left( \frac{2p-1}{p}\right) q_{1}}\left\vert
E_{k,s}^{1}\right\vert ^{\frac{q_{1}}{p}}\left(
\int\limits_{E_{k,s}^{1}}\left\vert \nabla \varphi \right\vert
^{p}\,dx\right) ^{\frac{\alpha }{p}}\right] $ \\ 
$\leq 2^{\left( \frac{p-1}{p}\right) q_{1}}\left( C_{SI}^{\alpha }\left\vert
E_{k,s}^{1}\right\vert ^{\Theta _{1}}\left\Vert a\right\Vert _{L^{\sigma
_{1}}\left( E_{k,s}^{1}\right) }\right) $ \\ 
$\cdot \left[ \frac{1}{\varepsilon ^{\frac{p+q_{1}}{p-\left( \alpha
+q_{1}\right) }}}+\varepsilon \int\limits_{E_{k,s}^{1}}\left\vert \nabla
\varphi \right\vert ^{p}\,dx+\frac{2^{\left( \frac{2p-1}{p}\right) q_{1}}}{%
\varepsilon ^{\frac{\alpha }{p-\alpha }}}\left(
\int\limits_{E_{k,s}^{1}}\left\vert \nabla u\right\vert ^{p}\,dx\right)
^{\left( \frac{q_{1}}{p}\right) \frac{p}{p-\alpha }}+\varepsilon 2^{\left( 
\frac{2p-1}{p}\right) q_{1}}\int\limits_{E_{k,s}^{1}}\left\vert \nabla
\varphi \right\vert ^{p}\,dx\right. $ \\ 
$\left. +\frac{2^{\left( \frac{2p-1}{p}\right) q_{1}}}{\varepsilon ^{\frac{%
\alpha }{p-\alpha }}}\left\vert E_{k,s}^{1}\right\vert ^{\left( \frac{q_{1}}{%
p}\right) \frac{p}{p-\alpha }}+\varepsilon 2^{\left( \frac{2p-1}{p}\right)
q_{1}}\int\limits_{E_{k,s}^{1}}\left\vert \nabla \varphi \right\vert ^{p}\,dx%
\right] $ \\ 
$\leq \varepsilon \left( 1+2^{1+\left( \frac{2p-1}{p}\right) q_{1}}\right)
2^{\left( \frac{p-1}{p}\right) q_{1}}\left( C_{SI}^{\alpha }\left\vert
E_{k,s}^{1}\right\vert ^{\Theta _{1}}\left\Vert a\right\Vert _{L^{\sigma
_{1}}\left( E_{k,s}^{1}\right) }\right) \int\limits_{E_{k,s}^{1}}\left\vert
\nabla \varphi \right\vert ^{p}\,dx$ \\ 
$+2^{\left( \frac{p-1}{p}\right) q_{1}}\left[ \frac{1}{\varepsilon ^{\frac{%
p+q_{1}}{p-\left( \alpha +q_{1}\right) }}}+\frac{2^{\left( \frac{2p-1}{p}%
\right) q_{1}}}{\varepsilon ^{\frac{\alpha }{p-\alpha }}}\left(
\int\limits_{E_{k,s}^{1}}\left\vert \nabla u\right\vert ^{p}\,dx\right)
^{\left( \frac{q_{1}}{p}\right) \frac{p}{p-\alpha }}+\frac{2^{\left( \frac{%
2p-1}{p}\right) q_{1}}}{\varepsilon ^{\frac{\alpha }{p-\alpha }}}\left\vert
E_{k,s}^{1}\right\vert ^{\left( \frac{q_{1}}{p}\right) \frac{p}{p-\alpha }}%
\right] $ \\ 
$\cdot \left( C_{SI}^{\alpha }\left\vert E_{k,s}^{1}\right\vert ^{\Theta
_{1}}\left\Vert a\right\Vert _{L^{\sigma _{1}}\left( E_{k,s}^{1}\right)
}\right) $ \\ 
$\leq \varepsilon \left( 1+2^{1+\left( \frac{2p-1}{p}\right) q_{1}}\right)
2^{\left( \frac{p-1}{p}\right) q_{1}}\left( C_{SI}^{\alpha }\left\vert
\Sigma \right\vert ^{\Theta _{1}}\left\Vert a\right\Vert _{L^{\sigma
_{1}}\left( \Sigma \right) }\right) \int\limits_{E_{k,s}^{1}}\left\vert
\nabla \varphi \right\vert ^{p}\,dx$ \\ 
$+2^{\left( \frac{p-1}{p}\right) q_{1}}\left[ \frac{1}{\varepsilon ^{\frac{%
p+q_{1}}{p-\left( \alpha +q_{1}\right) }}}+\frac{2^{\left( \frac{2p-1}{p}%
\right) q_{1}}}{\varepsilon ^{\frac{\alpha }{p-\alpha }}}\left(
\int\limits_{\Sigma }\left\vert \nabla u\right\vert ^{p}\,dx\right) ^{\left( 
\frac{q_{1}}{p}\right) \frac{p}{p-\alpha }}+\frac{2^{\left( \frac{2p-1}{p}%
\right) q_{1}}}{\varepsilon ^{\frac{\alpha }{p-\alpha }}}\left\vert \Sigma
\right\vert ^{\left( \frac{q_{1}}{p}\right) \frac{p}{p-\alpha }}\right] $ \\ 
$\cdot \left( C_{SI}^{\alpha }\left\vert E_{k,s}^{1}\right\vert ^{\Theta
_{1}}\left\Vert a\right\Vert _{L^{\sigma _{1}}\left( E_{k,s}^{1}\right)
}\right) $%
\end{tabular}%
\end{equation}%
and%
\begin{equation*}
\begin{tabular}{l}
$\left( \int\limits_{E_{k,s}^{1}}\left\vert \nabla \varphi \right\vert
^{p}\,dx\right) ^{\frac{\beta }{p}}\left( \int\limits_{E_{k,s}^{1}}\left(
\left\vert \varphi \right\vert +2\left\vert u\right\vert +1\right) ^{p^{\ast
}}\,dx\right) ^{\frac{q_{2}}{p^{\ast }}}\left( \left\vert
E_{k,s}^{1}\right\vert ^{\Theta _{2}}\left\Vert b\right\Vert _{L^{\sigma
_{2}}\left( E_{k,s}^{1}\right) }\right) $ \\ 
$\leq \left( \int\limits_{E_{k,s}^{1}}\left\vert \nabla \varphi \right\vert
^{p}\,dx\right) ^{\frac{\beta }{p}}\left( \left\vert E_{k,s}^{1}\right\vert
^{\Theta _{2}}\left\Vert b\right\Vert _{L^{\sigma _{2}}\left(
E_{k,s}^{1}\right) }\right) $ \\ 
$\cdot \left( 2^{p^{\ast }-1}\int\limits_{E_{k,s}^{1}}\left\vert \varphi
\right\vert ^{p^{\ast }}\,dx+2^{3p^{\ast
}-2}\int\limits_{E_{k,s}^{1}}\left\vert u\right\vert ^{p^{\ast
}}\,dx+2^{3p^{\ast }-2}\left\vert E_{k,s}^{1}\right\vert \right) ^{\frac{%
q_{2}}{p^{\ast }}}$ \\ 
$\leq \left( \int\limits_{E_{k,s}^{1}}\left\vert \nabla \varphi \right\vert
^{p}\,dx\right) ^{\frac{\beta }{p}}\left( \left\vert E_{k,s}^{1}\right\vert
^{\Theta _{2}}\left\Vert b\right\Vert _{L^{\sigma _{2}}\left(
E_{k,s}^{1}\right) }\right) $ \\ 
$\cdot \left[ 2^{\left( p^{\ast }-1\right) \frac{q_{2}}{p^{\ast }}}\left(
\int\limits_{E_{k,s}^{1}}\left\vert \varphi \right\vert ^{p^{\ast
}}\,dx\right) ^{\frac{q_{2}}{p^{\ast }}}+2^{\left( 3p^{\ast }-2\right) \frac{%
q_{2}}{p^{\ast }}}\left( \int\limits_{E_{k,s}^{1}}\left\vert u\right\vert
^{p^{\ast }}\,dx\right) ^{\frac{q_{2}}{p^{\ast }}}+2^{\left( 3p^{\ast
}-2\right) \frac{q_{2}}{p^{\ast }}}\left\vert E_{k,s}^{1}\right\vert ^{\frac{%
q_{2}}{p^{\ast }}}\right] $ \\ 
$\leq \left( \int\limits_{E_{k,s}^{1}}\left\vert \nabla \varphi \right\vert
^{p}\,dx\right) ^{\frac{\beta }{p}}\left( \left\vert E_{k,s}^{1}\right\vert
^{\Theta _{2}}\left\Vert b\right\Vert _{L^{\sigma _{2}}\left(
E_{k,s}^{1}\right) }\right) $ \\ 
$\cdot \left[ 2^{\left( p^{\ast }-1\right) \frac{q_{2}}{p^{\ast }}%
}C_{SI}^{q_{2}}\left( \int\limits_{E_{k,s}^{1}}\left\vert \nabla \varphi
\right\vert ^{p}\,dx\right) ^{\frac{q_{2}}{p}}+2^{\left( 3p^{\ast }-2\right) 
\frac{q_{2}}{p^{\ast }}}\left( \int\limits_{E_{k,s}^{1}}\left\vert
u\right\vert ^{p^{\ast }}\,dx\right) ^{\frac{q_{2}}{p^{\ast }}}+2^{\left(
3p^{\ast }-2\right) \frac{q_{2}}{p^{\ast }}}\left\vert
E_{k,s}^{1}\right\vert ^{\frac{q_{2}}{p^{\ast }}}\right] $ \\ 
$\leq \left( \left\vert E_{k,s}^{1}\right\vert ^{\Theta _{2}}\left\Vert
b\right\Vert _{L^{\sigma _{2}}\left( E_{k,s}^{1}\right) }\right) $ \\ 
$\cdot \left[ 2^{\left( p^{\ast }-1\right) \frac{q_{2}}{p^{\ast }}%
}C_{SI}^{q_{2}}\left( \int\limits_{E_{k,s}^{1}}\left\vert \nabla \varphi
\right\vert ^{p}\,dx\right) ^{\frac{\beta +q_{2}}{p}}+2^{\left( 3p^{\ast
}-2\right) \frac{q_{2}}{p^{\ast }}}\left(
\int\limits_{E_{k,s}^{1}}\left\vert u\right\vert ^{p^{\ast }}\,dx\right) ^{%
\frac{q_{2}}{p^{\ast }}}\left( \int\limits_{E_{k,s}^{1}}\left\vert \nabla
\varphi \right\vert ^{p}\,dx\right) ^{\frac{\beta }{p}}\right. $ \\ 
$\left. +2^{\left( 3p^{\ast }-2\right) \frac{q_{2}}{p^{\ast }}}\left\vert
E_{k,s}^{1}\right\vert ^{\frac{q_{2}}{p^{\ast }}}\left(
\int\limits_{E_{k,s}^{1}}\left\vert \nabla \varphi \right\vert
^{p}\,dx\right) ^{\frac{\beta }{p}}\right] $ \\ 
$\leq \left( \left\vert E_{k,s}^{1}\right\vert ^{\Theta _{2}}\left\Vert
b\right\Vert _{L^{\sigma _{2}}\left( E_{k,s}^{1}\right) }\right) $ \\ 
$\cdot \left[ \frac{2^{\left( p^{\ast }-1\right) \frac{q_{2}}{p^{\ast }}}}{%
\varepsilon ^{\frac{\beta +q_{2}}{p-\left( \beta +q_{2}\right) }}}%
C_{SI}^{q_{2}}+\varepsilon 2^{\left( p^{\ast }-1\right) \frac{q_{2}}{p^{\ast
}}}C_{SI}^{q_{2}}\int\limits_{E_{k,s}^{1}}\left\vert \nabla \varphi
\right\vert ^{p}\,dx\right. $ \\ 
$\left. +\frac{2^{\left( 3p^{\ast }-2\right) \frac{q_{2}}{p^{\ast }}}}{%
\varepsilon ^{\frac{\beta }{p-\beta }}}\left(
\int\limits_{E_{k,s}^{1}}\left\vert u\right\vert ^{p^{\ast }}\,dx\right)
^{\left( \frac{q_{2}}{p^{\ast }}\right) \frac{p}{p-\beta }}+\varepsilon
2^{\left( 3p^{\ast }-2\right) \frac{q_{2}}{p^{\ast }}}\int%
\limits_{E_{k,s}^{1}}\left\vert \nabla \varphi \right\vert ^{p}\,dx\right. $
\\ 
$\left. +\frac{2^{\left( 3p^{\ast }-2\right) \frac{q_{2}}{p^{\ast }}}}{%
\varepsilon ^{\frac{\beta }{p-\beta }}}\left\vert E_{k,s}^{1}\right\vert
^{\left( \frac{q_{2}}{p^{\ast }}\right) \frac{p}{p-\beta }}+\varepsilon
2^{\left( 3p^{\ast }-2\right) \frac{q_{2}}{p^{\ast }}}\int%
\limits_{E_{k,s}^{1}}\left\vert \nabla \varphi \right\vert ^{p}\,dx\right] $
\\ 
$\leq \varepsilon \left( \left\vert \Sigma \right\vert ^{\Theta
_{2}}\left\Vert b\right\Vert _{L^{\sigma _{2}}\left( \Sigma \right) }\right)
\left( 2^{\left( p^{\ast }-1\right) \frac{q_{2}}{p^{\ast }}%
}C_{SI}^{q_{2}}+2^{1+\left( 3p^{\ast }-2\right) \frac{q_{2}}{p^{\ast }}%
}\right) \int\limits_{E_{k,s}^{1}}\left\vert \nabla \varphi \right\vert
^{p}\,dx$ \\ 
$+\left( \left\vert E_{k,s}^{1}\right\vert ^{\Theta _{2}}\left\Vert
b\right\Vert _{L^{\sigma _{2}}\left( \Sigma \right) }\right) \left[ \frac{%
2^{\left( p^{\ast }-1\right) \frac{q_{2}}{p^{\ast }}}}{\varepsilon ^{\frac{%
\beta +q_{2}}{p-\left( \beta +q_{2}\right) }}}C_{SI}^{q_{2}}+\frac{2^{\left(
3p^{\ast }-2\right) \frac{q_{2}}{p^{\ast }}}}{\varepsilon ^{\frac{\beta }{%
p-\beta }}}\left( \int\limits_{B_{s}\left( x_{0}\right) }\left\vert
u\right\vert ^{p^{\ast }}\,dx\right) ^{\left( \frac{q_{2}}{p^{\ast }}\right) 
\frac{p}{p-\beta }}+\frac{2^{\left( 3p^{\ast }-2\right) \frac{q_{2}}{p^{\ast
}}}}{\varepsilon ^{\frac{\beta }{p-\beta }}}\left\vert \Sigma \right\vert
^{\left( \frac{q_{2}}{p^{\ast }}\right) \frac{p}{p-\beta }}\right] $%
\end{tabular}%
\end{equation*}%
then using the Embendding Sobolev Theorem we get%
\begin{equation*}
\begin{tabular}{l}
$\int\limits_{E_{k,s}^{1}}a\left( x\right) \left\vert \varphi \right\vert
^{\alpha }\left( \left\vert \nabla u+\nabla \varphi \right\vert +\left\vert
\nabla u\right\vert +1\right) ^{q_{1}}+b\left( x\right) \left( \left\vert
u+\varphi \right\vert +\left\vert u\right\vert +1\right) ^{q_{2}}\left\vert
\nabla \varphi \right\vert ^{\beta }\,dx$ \\ 
$\leq \varepsilon D_{1,\Sigma }\int\limits_{E_{k,s}^{1}}\left\vert \nabla
\varphi \right\vert ^{p}\,dx+D_{2,\varepsilon ,\Sigma }\left\vert
E_{k,s}^{1}\right\vert ^{\Theta _{1}}+D_{3,\varepsilon ,\Sigma }\left\vert
E_{k,s}^{1}\right\vert ^{\Theta _{2}}$%
\end{tabular}%
\end{equation*}%
where%
\begin{eqnarray*}
D_{1,\Sigma } &=&\left( 1+2^{1+\left( \frac{2p-1}{p}\right) q_{1}}\right)
2^{\left( \frac{p-1}{p}\right) q_{1}}\left( C_{SI}^{\alpha }\left\vert
\Sigma \right\vert ^{\Theta _{1}}\left\Vert a\right\Vert _{L^{\sigma
_{1}}\left( \Sigma \right) }\right) \\
&&+\left( \left\vert \Sigma \right\vert ^{\Theta _{2}}\left\Vert
b\right\Vert _{L^{\sigma _{2}}\left( \Sigma \right) }\right) \left(
2^{\left( p^{\ast }-1\right) \frac{q_{2}}{p^{\ast }}}C_{SI}^{q_{2}}+2^{1+%
\left( 3p^{\ast }-2\right) \frac{q_{2}}{p^{\ast }}}\right)
\end{eqnarray*}%
\begin{equation*}
D_{2,\varepsilon ,\Sigma }=C_{SI}^{\alpha }\left\Vert a\right\Vert
_{L^{\sigma _{1}}\left( \Sigma \right) }\left[ \frac{1}{\varepsilon ^{\frac{%
p+q_{1}}{p-\left( \alpha +q_{1}\right) }}}+\frac{2^{\left( \frac{2p-1}{p}%
\right) q_{1}}}{\varepsilon ^{\frac{\alpha }{p-\alpha }}}\left( \left\Vert
u\right\Vert _{W^{1,p}\left( \Sigma \right) }\right) ^{\frac{pq_{1}}{%
p-\alpha }}+\frac{2^{\left( \frac{2p-1}{p}\right) q_{1}}}{\varepsilon ^{%
\frac{\alpha }{p-\alpha }}}\left\vert \Sigma \right\vert ^{\left( \frac{q_{1}%
}{p}\right) \frac{p}{p-\alpha }}\right]
\end{equation*}%
and%
\begin{equation*}
D_{3,\varepsilon ,\Sigma }=\left\Vert b\right\Vert _{L^{\sigma _{2}}\left(
\Sigma \right) }\left[ \frac{2^{\left( p^{\ast }-1\right) \frac{q_{2}}{%
p^{\ast }}}}{\varepsilon ^{\frac{\beta +q_{2}}{p-\left( \beta +q_{2}\right) }%
}}C_{SI}^{q_{2}}+\frac{2^{\left( 3p^{\ast }-2\right) \frac{q_{2}}{p^{\ast }}}%
}{\varepsilon ^{\frac{\beta }{p-\beta }}}C_{SE}^{\frac{pq_{2}}{p-\beta }%
}\left( \left\Vert u\right\Vert _{W^{1,p}\left( \Sigma \right) }\right) ^{%
\frac{pq_{2}}{p-\beta }}+\frac{2^{\left( 3p^{\ast }-2\right) \frac{q_{2}}{%
p^{\ast }}}}{\varepsilon ^{\frac{\beta }{p-\beta }}}\left\vert \Sigma
\right\vert ^{\left( \frac{q_{2}}{p^{\ast }}\right) \frac{p}{p-\beta }}%
\right]
\end{equation*}%
Since%
\begin{equation*}
\int\limits_{E_{k,s}^{1}}\left\vert \nabla \varphi \right\vert ^{p}\,dx\leq
\int\limits_{E_{k,s}^{1}}\eta ^{p}\left\vert \nabla u^{1}\right\vert
^{p}\,dx+2^{p}\int\limits_{E_{k,s}^{1}}\frac{\left( u^{1}-k\right) ^{p}}{%
\left( s-t\right) ^{p}}\,dx
\end{equation*}%
it follows%
\begin{equation}
\begin{tabular}{l}
$\int\limits_{E_{k,s}^{1}}a\left( x\right) \left\vert \varphi \right\vert
^{\alpha }\left( \left\vert \nabla u+\nabla \varphi \right\vert +\left\vert
\nabla u\right\vert +1\right) ^{q_{1}}+b\left( x\right) \left( \left\vert
u+\varphi \right\vert +\left\vert u\right\vert +1\right) ^{q_{2}}\left\vert
\nabla \varphi \right\vert ^{\beta }\,dx$ \\ 
$\leq \varepsilon D_{1,\Sigma }\int\limits_{E_{k,s}^{1}}\eta ^{p}\left\vert
\nabla u^{1}\right\vert ^{p}\,dx+\varepsilon 2^{p}D_{1,\Sigma
}\int\limits_{E_{k,s}^{1}}\frac{\left( u^{1}-k\right) ^{p}}{\left(
s-t\right) ^{p}}\,dx+D_{2,\varepsilon ,\Sigma }\left\vert
E_{k,s}^{1}\right\vert ^{\Theta _{1}}+D_{3,\varepsilon ,\Sigma }\left\vert
E_{k,s}^{1}\right\vert ^{\Theta _{2}}$%
\end{tabular}%
\end{equation}

\bigskip Usinig (3.12) and (3.16) we get%
\begin{equation*}
\begin{tabular}{l}
$\int\limits_{E_{k,s}^{1}}\left\vert \nabla u^{1}\right\vert ^{p}\,dx$ \\ 
$\leq 2^{p-1}\int\limits_{E_{k,s}^{1}}\left( 1-\eta ^{p}\right) \left\vert
\nabla u^{1}\right\vert
^{p}\,dx+2^{2p-1}p^{p}\int\limits_{E_{k,s}^{1}-E_{k,t}^{1}}\frac{\left(
u^{1}-k\right) ^{p}}{\left( s-t\right) ^{p}}\,dx$ \\ 
$+\varepsilon D_{1,\Sigma }\int\limits_{E_{k,s}^{1}}\eta ^{p}\left\vert
\nabla u^{1}\right\vert ^{p}\,dx+\varepsilon 2^{p}D_{1,\Sigma
}\int\limits_{E_{k,s}^{1}}\frac{\left( u^{1}-k\right) ^{p}}{\left(
s-t\right) ^{p}}\,dx+D_{2,\varepsilon ,\Sigma }\left\vert
E_{k,s}^{1}\right\vert ^{\Theta _{1}}+D_{3,\varepsilon ,\Sigma }\left\vert
E_{k,s}^{1}\right\vert ^{\Theta _{2}}$%
\end{tabular}%
\end{equation*}

Fix $\varepsilon =\frac{1}{2D_{1,\Sigma }}$ \ it follows%
\begin{equation}
\begin{tabular}{l}
$\frac{1}{2}\int\limits_{E_{k,s}^{1}}\left\vert \nabla u^{1}\right\vert
^{p}\,dx$ \\ 
$\leq 2^{p-1}\int\limits_{E_{k,s}^{1}}\left( 1-\eta ^{p}\right) \left\vert
\nabla u^{1}\right\vert ^{p}\,dx+\left( 2^{p-1}+2^{2p-1}p^{p}\right)
\int\limits_{E_{k,s}^{1}-E_{k,t}^{1}}\frac{\left( u^{1}-k\right) ^{p}}{%
\left( s-t\right) ^{p}}\,dx$ \\ 
$+D_{2,\Sigma }\left\vert E_{k,s}^{1}\right\vert ^{\Theta _{1}}+D_{3,\Sigma
}\left\vert E_{k,s}^{1}\right\vert ^{\Theta _{2}}$%
\end{tabular}
\label{81}
\end{equation}%
where%
\begin{equation*}
D_{2,\Sigma }=C_{SI}^{\alpha }\left\Vert a\right\Vert _{L^{\sigma
_{1}}\left( \Sigma \right) }\left[ \frac{1}{\left( \frac{1}{2D_{1,\Sigma }}%
\right) ^{\frac{p+q_{1}}{p-\left( \alpha +q_{1}\right) }}}+\frac{2^{\left( 
\frac{2p-1}{p}\right) q_{1}}}{\left( \frac{1}{2D_{1,\Sigma }}\right) ^{\frac{%
\alpha }{p-\alpha }}}\left( \left\Vert u\right\Vert _{W^{1,p}\left( \Sigma
\right) }\right) ^{\frac{pq_{1}}{p-\alpha }}+\frac{2^{\left( \frac{2p-1}{p}%
\right) q_{1}}}{\left( \frac{1}{2D_{1,\Sigma }}\right) ^{\frac{\alpha }{%
p-\alpha }}}\left\vert \Sigma \right\vert ^{\left( \frac{q_{1}}{p}\right) 
\frac{p}{p-\alpha }}\right]
\end{equation*}%
and%
\begin{equation*}
D_{3,\Sigma }=\left\Vert b\right\Vert _{L^{\sigma _{2}}\left( \Sigma \right)
}\left[ \frac{2^{\left( p^{\ast }-1\right) \frac{q_{2}}{p^{\ast }}}}{\left( 
\frac{1}{2D_{1,\Sigma }}\right) ^{\frac{\beta +q_{2}}{p-\left( \beta
+q_{2}\right) }}}C_{SI}^{q_{2}}+\frac{2^{\left( 3p^{\ast }-2\right) \frac{%
q_{2}}{p^{\ast }}}}{\left( \frac{1}{2D_{1,\Sigma }}\right) ^{\frac{\beta }{%
p-\beta }}}C_{SE}^{\frac{pq_{2}}{p-\beta }}\left( \left\Vert u\right\Vert
_{W^{1,p}\left( \Sigma \right) }\right) ^{\frac{pq_{2}}{p-\beta }}+\frac{%
2^{\left( 3p^{\ast }-2\right) \frac{q_{2}}{p^{\ast }}}}{\left( \frac{1}{%
2D_{1,\Sigma }}\right) ^{\frac{\beta }{p-\beta }}}\left\vert \Sigma
\right\vert ^{\left( \frac{q_{2}}{p^{\ast }}\right) \frac{p}{p-\beta }}%
\right] \text{.}
\end{equation*}

\bigskip It is easy to observe that the constants $D_{2,\Sigma }$ and $%
D_{3,\Sigma }$ are independent of the point $x_{0}$ and that they depend
only on the initial data. \ 

Using (\ref{81}), we get%
\begin{equation}
\begin{tabular}{l}
$\int\limits_{E_{k,s}^{1}}\left\vert \nabla u^{1}\right\vert ^{p}\,dx$ \\ 
$\leq \frac{2^{p}}{1+2^{p}}\int\limits_{E_{k,s}^{1}}\left\vert \nabla
u^{1}\right\vert ^{p}\,dx+\left( \frac{2^{p}+2^{2p}p^{p}}{1+2^{p}}\right)
\int\limits_{E_{k,s}^{1}-E_{k,t}^{1}}\frac{\left( u^{1}-k\right) ^{p}}{%
\left( s-t\right) ^{p}}\,dx$ \\ 
$+\frac{2D_{2,\Sigma }}{1+2^{p}}\left\vert E_{k,s}^{1}\right\vert ^{\Theta
_{1}}+\frac{2D_{3,\Sigma }}{1+2^{p}}\left\vert E_{k,s}^{1}\right\vert
^{\Theta _{2}}$%
\end{tabular}
\label{83}
\end{equation}

\bigskip Now, using Lemma \ref{lem3} we get%
\begin{equation}
\begin{tabular}{l}
$\int\limits_{A_{k,\varrho }^{1}}\left\vert \nabla u^{1}\right\vert
^{p}\,dx\leq \frac{C_{C,1}}{\left( R-\varrho \right) ^{p}}%
\int\limits_{A_{k,R}^{1}}\left( u^{1}-k\right) ^{p}\,dx+C_{C,2}\left\vert
A_{k,s}^{1}\right\vert ^{\Theta _{1}}+C_{C,3}\left\vert
A_{k,s}^{1}\right\vert ^{\Theta _{2}}$%
\end{tabular}%
\end{equation}

Since $-u$ is a local minimizer of the following integral functional

\bigskip 
\begin{equation*}
\tilde{J}\left( v,\Omega \right) =\int\limits_{\Omega }\sum\limits_{\alpha
=1}^{n}\left\vert \nabla v^{\alpha }\left( x\right) \right\vert ^{p}+\tilde{G%
}\left( x,v\left( x\right) ,\nabla v\left( x\right) \right) \,dx
\end{equation*}%
where $\tilde{G}\left( x,v\left( x\right) ,\nabla v\left( x\right) \right)
=G\left( x,-v\left( x\right) ,-\nabla v\left( x\right) \right) $ then we get%
\begin{equation}
\begin{tabular}{l}
$\int\limits_{B_{k,\varrho }^{1}}\left\vert \nabla u^{1}\right\vert
^{p}\,dx\leq \frac{C_{C,1}}{\left( R-\varrho \right) ^{p}}%
\int\limits_{B_{k,R}^{1}}\left( k-u^{1}\right) ^{p}\,dx+C_{C,2}\left\vert
B_{k,s}^{1}\right\vert ^{\Theta _{1}}+C_{C,3}\left\vert
B_{k,s}^{1}\right\vert ^{\Theta _{2}}$%
\end{tabular}%
\end{equation}%
Similarly we can proceed for $u^{\alpha }$ with $\alpha =2,...,m$. Since $%
\Theta _{1},\Theta _{2}>1-\frac{p}{n}+\epsilon $ then it follows%
\begin{equation*}
\int\limits_{A_{k,\varrho }^{1}}\left\vert \nabla u^{\alpha }\right\vert
^{p}\,dx\leq \frac{C_{C,1}}{\left( R-\varrho \right) ^{p}}%
\int\limits_{A_{k,R}^{1}}\left( u^{\alpha }-k\right) ^{p}\,dx+\left(
C_{C,2}+C_{C,3}\right) \left\vert A_{k,s}^{1}\right\vert ^{1-\frac{p}{n}%
+\epsilon }
\end{equation*}%
and%
\begin{equation*}
\int\limits_{B_{k,\varrho }^{1}}\left\vert \nabla u^{\alpha }\right\vert
^{p}\,dx\leq \frac{C_{C,1}}{\left( R-\varrho \right) ^{p}}%
\int\limits_{B_{k,R}^{1}}\left( k-u^{\alpha }\right) ^{p}\,dx+\left(
C_{C,2}+C_{C,3}\right) \left\vert B_{k,s}^{1}\right\vert ^{1-\frac{p}{n}%
+\epsilon }
\end{equation*}%
for every $\alpha =1,...,m$.

\bigskip

\section{Proof of Theorem 1}

\bigskip The proof follows by applying Theorem 2 and Theorem 8.

\end{document}